\pdfoutput=1
\documentclass[a4paper,11pt]{article}
\usepackage[english]{babel}
\usepackage{amsfonts}
\usepackage{amsthm}
\usepackage{amsmath}
\usepackage{amssymb}
\usepackage{enumerate}
\usepackage{mathtools}
\usepackage{tikz}
\usepackage{tikz-cd}
\usepackage[all]{xy}
\usepackage{graphicx}
\usepackage{amsmath}
\usepackage{adjustbox}
\usepackage{graphicx}
\usepackage{extpfeil}
\usepackage{comment}
\usepackage{float}
\usepackage[labelformat=empty]{caption}
\usepackage[toc]{appendix}
\usepackage[left=3.2cm,top=3cm,right=3.2cm,bottom=3cm]{geometry}
\usepackage{resizegather}
\usepackage[activate={true, nocompatibility}, final, tracking=true, kerning=true, spacing=true, factor=1100, stretch=10, shrink=10]{microtype}
\usepackage{authblk}
\usepackage{tikz}
\usepackage{dsfont}
\usetikzlibrary{matrix}

\usepackage{color}
\usepackage{resizegather}
\usepackage[normalem]{ulem}
\usepackage{enotez}
\setenotez{backref=true}

\usepackage{hyperref}
\hypersetup{
colorlinks=true,
linkcolor=[RGB]{0,0,204},
filecolor=magenta,      
urlcolor=cyan,
citecolor=[RGB]{0,204,0},
}
\usepackage[nameinlink]{cleveref}

\usepackage{mathrsfs}
\usepackage{mathscinet}

\theoremstyle{plain}
\newtheorem{thm}{Theorem}[section]
\newtheorem{coroll}[thm]{Corollary}
\newtheorem{defn}[thm]{Definition}
\newtheorem{lemma}[thm]{Lemma}

\newtheorem{prop}[thm]{Proposition}

\newtheorem{notn}[thm]{Notation}
\newtheorem*{thm*}{Main Theorem}
\newtheorem*{defn*}{Definition}

\theoremstyle{definition}
\newtheorem{example}[thm]{Example}
\newtheorem{remark}[thm]{Remark}

\newtheorem{introthm}{Theorem}

\newcommand{\dash}{{\operatorname{-}}}

\tikzset{
  symbol/.style={
    draw=none,
    every to/.append style={
      edge node={node [sloped, allow upside down, auto=false]{$#1$}}}
  }
}
\microtypecontext{spacing=nonfrench}

\makeatletter % changes the catcode of @ to 11
\def\makeCal#1{%
\expandafter\newcommand\csname c#1\endcsname{\mathcal{#1}}}

\def\makeBB#1{%
\expandafter\newcommand\csname b#1\endcsname{\mathbb{#1}}}

\def\makeFrak#1{%
\expandafter\newcommand\csname f#1\endcsname{\mathfrak{#1}}}

\def\makeScr#1{%
\expandafter\newcommand\csname s#1\endcsname{\mathscr{#1}}}

\count@=0
\loop
\advance\count@ 1
\edef\y{\@Alph\count@} % Alph is upper case 
\expandafter\makeCal\y
\expandafter\makeBB\y
\expandafter\makeFrak\y
\expandafter\makeScr\y
\ifnum\count@<26
\repeat
\def\makelowercaseFrak#1{%
\expandafter\newcommand\csname mf#1\endcsname{\mathfrak{#1}}}
\count@=0
\loop
\advance\count@ 1
\edef\y{\@alph\count@} % alph is lower case.. 
\expandafter\makelowercaseFrak\y
\ifnum\count@<26
\repeat
\DeclareMathOperator{\GL}{GL}
\DeclareMathOperator{\Bun}{Bun}
\DeclareMathOperator{\ev}{ev}

\DeclareMathOperator{\wt}{wt}

\DeclareMathOperator{\Sym}{Sym}

\DeclareMathOperator{\Filt}{Filt}

\DeclareMathOperator{\Map}{\operatorname{Map}}

\DeclareMathOperator{\Gr}{Gr}

\DeclareMathOperator{\Tot}{Tot}
\DeclareMathOperator{\stab}{stab}
\DeclareMathOperator{\Sch}{Sch}
\DeclareMathOperator{\tr}{tr}
\DeclareMathOperator{\rk}{rk}

\newcommand{\F}{\mathcal{F}}

\newcommand{\on}{\operatorname}

\newcommand{\Coh}{\mathscr{C}\!oh}

\newcommand{\Aut}{ \on{Aut} } 

\newcommand{\Hom}{ \on{Hom}}

\newcommand{\Spec}{\on{Spec}}

\newcommand{\Pic}{\on{Pic}}

\newcommand\dirac{/\kern-1.2ex\partial} % Dirac operator
\newcommand\qu{/\kern-.7ex/} % Categorical quotients
\newcommand\lqu{\backslash \kern-.7ex \backslash} % Categorical
\newcommand\dr{r_+ \kern-.7ex - \kern-.7ex r_-}
 
\tikzset{
  symbol/.style={
    draw=none,
    every to/.append style={
      edge node={node [sloped, allow upside down, auto=false]{$#1$}}}
  }
}

    \newtheoremstyle{TheoremNum}
        {\topsep}{\topsep}              %%% space between body and thm
        {\itshape}                      %%% Thm body font
        {}                              %%% Indent amount (empty = no indent)
        {\bfseries}                     %%% Thm head font
        {.}                             %%% Punctuation after thm head
        { }                             %%% Space after thm head
        {\thmname{#1}\thmnote{ \bfseries #3}}%%% Thm head spec
    \theoremstyle{TheoremNum}

\newcommand{\GBun}{\mathscr{G}\!Bun}
\newcommand{\Higgs}{\mathscr{H}\!iggs}
\newcommand{\GHiggs}{\mathscr{G}\!\mathscr{H}\!iggs}
\newcommand{\MHiggs}{Higgs}
\newcommand{\MGHiggs}{GHiggs}

\begin{document}
\title{\textbf{The meromorphic Hitchin fibration over stable pointed curves: moduli spaces}}
\author{Ron Donagi and Andres Fernandez Herrero}
\date{}
\maketitle
\begin{abstract}
We construct a universal partial compactification of the relative moduli space of semistable meromorphic Higgs bundles over the stack of stable pointed curves. It parametrizes meromorphic Gieseker Higgs bundles, and is equipped with a flat and proper extension of the usual Hitchin morphism. Over an open subset of the Hitchin base parametrizing allowable nodal spectral covers, we describe the relation of the fibers to compactified Jacobians, thus establishing an analogue of the BNR correspondence. We also construct a version of the moduli space where we require the residues of the meromorphic Higgs bundle to lie in a given set of nilpotent conjugacy classes. In this latter case, we show that there is a flat and proper Hitchin morphism to the flat degeneration of the corresponding family of Hitchin bases constructed in previous physics work.
\end{abstract}
\tableofcontents

\begin{section}{Introduction}

Given a smooth projective connected curve $C$ over a fixed ground field $k$ of characteristic $0$, the Dolbeault moduli space $\MHiggs_{C}^{N,\chi}$ \cite{hitchin_self_duality,simpson-repnII}
parametrizes semistable Higgs bundles $(E, \phi: E \to E \otimes \omega_C)$, where $E$ is a vector bundle on $C$ of some given fixed rank $N$ and fixed Euler characteristic $\chi (E) = \chi$.
Hitchin’s morphism:
\[
H: \MHiggs_C^{N,\chi} \to A_C
\]
sends a Higgs bundle to the characteristic polynomial of its Higgs field $\phi$. Here $A_C$ is a vector space parametrizing the coefficients of these characteristic polynomials,
\[
A_C := \oplus_{i=1}^N H^0(C,\omega_C^{\otimes i}).
\]
The Hitchin morphism $H: \MHiggs_{C}^{N,\chi} \to A_C$ is an integrable system: there is a natural symplectic structure on $\MHiggs_C^{N,\chi}$, and $H$ is proper and Lagrangian.

There is a meromorphic version of the Dolbeault moduli space. In the meromorphic version, the curve $C$ is replaced by an $n$-pointed (aka: punctured) curve $(C,D)$, where $C$ is as before and  $D=\sum_{i=1}^n \sigma_i$ is a reduced divisor on $C$. Higgs bundles are replaced by meromorphic Higgs bundles $(E, \phi: E \to E \otimes \omega_C (D))$, where the Higgs field $\phi$ is now allowed to have (first order) poles at the punctures $\sigma_i$. Semistable meromorphic Higgs bundles are parametrized by a moduli space $\MHiggs_{C,D}^{N,\chi}$. The Hitchin base is replaced by
\[
A_{C,D} := \oplus_{i=1}^N H^0(C,(\omega_C (D))^{\otimes i}),
\]
and the Hitchin morphism is now $H: \MHiggs_{C,D}^N \to A_{C,D}$. In this meromorphic context, $\MHiggs_{C,D}^{N,\chi}$ has a natural Poisson structure, and the Hitchin morphism is still Lagrangian (in the Poisson sense) and proper, as shown in the work of Markman \cite{markman-94} and Bottacin \cite{bottacin-95}. 

The symplectic leaves in $\MHiggs_{C,D}^{N,\chi}$ are obtained by imposing that the residue of $\phi$ at each puncture $\sigma_i$ lives in a specified conjugacy class in $\mathfrak{gl}_N$. For each $n$-tuple $\cO_{\bullet}= (\cO_i)_{i=1}^n$ of nilpotent conjugacy classes, there is an integrable system 
\[H:\overline{\MHiggs}^{\cO_{\bullet}, \chi}_{C,D} \to A^{\cO_{\bullet}}_{C,D},\]
where $\overline{\MHiggs}^{\cO_{\bullet}, \chi}_{C,D}$ is the closure of the locus of semistable meromorphic Higgs bundles with residue classes given by $\cO_{\bullet}$, and $A^{\cO_{\bullet}}_{C,D}$ is a linear subspace of $A_{C,D}$.

All of these constructions work in the relative setting over families of smooth curves. If we denote by $\cM_{g,n}$ the Deligne-Mumford stack of stable smooth $n$-pointed curves, then there is a relative moduli space of semistable meromorphic Higgs bundles $\MHiggs^{N,\chi}_{g,n} \to \cM_{g,n}$, a vector bundle family of Hitchin bases $\cA^{\circ} \to {\cM}_{g,n}$,
and a (universal) Hitchin morphism $\MHiggs^{N,\chi}_{g,n} \to \cA^{\circ}$. 
Furthermore, given any fixed $n$-tuple $\cO_{\bullet}$ of nilpotent conjugacy classes, we have a family of symplectic systems 
\[ H:\overline{\MHiggs}^{\cO_{\bullet},\chi}_{g,n} \to \cA^{\cO_{\bullet},\circ},\]
where $\cA^{\cO_{\bullet},\circ} \to \cM_{g,n}$ is a vector subbundle of $\cA^{\circ} \to \cM_{g,n}$.

For various reasons, some coming from conformal field theory (see e.g. \cite{gaiotto-dualities} and \cite{BDD-modified-hitchin-gln}) it is natural to allow the (punctured) curve $(C,\sigma_{\bullet})$ to degenerate, i.e. to try to extend the family of Poisson systems and the corresponding symplectic systems for the symplectic leaves to the Deligne-Mumford compactification 
$\overline{\cM}_{g,n}$.

Whereas the vector bundle $\cA^{\circ}$ naively extends to a vector bundle $\cA$ over $\overline{\cM}_{g,n}$, the same is not true for the subbundle $\cA^{\cO_{\bullet},\circ}$. A first step in this direction was taken in \cite{BDD-modified-hitchin-gln}, where it was observed that a sequence of twists by certain boundary
divisors of $\overline{\cM}_{g,n}$ allows us to define a vector bundle
\[
\cA^{\cO_{\bullet}} \to \overline{\cM}_{g,n}
\]
whose restriction to the interior $\cM_{g,n}$ recovers the naïve family given by $\cA^{\cO_{\bullet},\circ}$.

Our first main result in this article is a construction of extensions of the different versions of the (meromorphic) Dolbeault moduli space to the compactification $\overline{\cM}_{g,n}$.
\begin{introthm}[=Theorems \ref{thm: theta-stratification and moduli space for logarithmic gieseker higgs} and \ref{thm: moduli space for closure of stack with fixed residues}] \label{introthm: moduli thm}
    Fix a ground field $k$ of characteristic $0$, a rank $N\geq 0$, and an Euler characteristic $\chi$. Then, the following hold:
    \begin{enumerate}[(1)]
        \item There is a moduli-theoretic flat extension $\MGHiggs_{g,n}^{N,\chi} \to \overline{\cM}_{g,n}$ of the relative Dolbeault moduli space $\MHiggs_{g,n}^{N,\chi} \to \cM_{g,n}$. There is an extension $H: \MGHiggs_{g,n}^{N,\chi} \to \cA$  of the (universal) Hitchin morphism such that $H$ is flat and proper.

        \item Given a fixed $n$-tuple $\cO_{\bullet}$ of nilpotent conjugacy classes, there is a moduli-theoretic flat extension $\overline{\MGHiggs}^{\cO_{\bullet}, \chi}_{g,n} \to \overline{\cM}_{g,n}$ of the family $\overline{\MHiggs}^{\cO_{\bullet},\chi}_{g,n} \to \cM_{g,n}$ of closed symplectic leaves. There is an extension  $H: \overline{\MGHiggs}^{\cO_{\bullet}, \chi}_{g,n} \to \cA^{\cO_{\bullet}}$ of the corresponding Hitchin morphism such that $H$ is flat and proper.
    \end{enumerate}
\end{introthm}

The meaning of ``moduli-theoretic" in the statement of \Cref{introthm: moduli thm} is that both $\MGHiggs_{g,n}^{N,\chi}$ and $\overline{\MGHiggs}^{\cO_{\bullet}, \chi}_{g,n}$ are relative moduli spaces of certain objects called meromorphic Gieseker Higgs bundles. 
\begin{defn}
    A meromorphic Gieseker Higgs bundle is a tuple $(C, \sigma_{\bullet}, E, \psi)$, where:
    \begin{enumerate}[(1)]
        \item $(C, \sigma_{\bullet})$ is a semistable $n$-pointed curve (i.e. an $n$-pointed prestable curve such that its pointed stabilization morphism $\varphi: (C, \sigma_{\bullet}) \to (\overline{C}, \overline{\sigma}_{\bullet})$ preserves the log-canonical bundles: $\varphi^*\left(\omega_{\overline{C}}(\sum_{i=1}^n \overline{\sigma}_i)\right) = \omega_{C}(\sum_{i=1}^n \sigma_i)$).

        \item $E$ is a vector bundle on $C$ satisfying the following:
\begin{enumerate}[(G1)]
\item The counit $\varphi^*\varphi_*(E) \to E$ of the stabilization morphism $\varphi$ is surjective.
\item The determinant $\det(E)$ is $\varphi$-ample.
\item The pushforward $\varphi_*(E)$ is a torsion-free sheaf on the stabilization $\overline{C}$.
\end{enumerate}

\item $\psi: E \to E \otimes \omega_{C}(\sum_{i=1}^n \sigma_i)$ is a meromorphic Higgs field. 
    \end{enumerate}
\end{defn}
The data of $(C, \sigma_{\bullet}, E)$ is usually known in the literature as a Gieseker vector bundle. We refer the reader to \Cref{remark: Gieseker bundle concrete} for a more explicit description of the conditions (G1), (G2) and (G3) above. The idea of using semistable modifications $(C, \sigma_{\bullet})$ in order to compactify the moduli of vector bundles on nodal curves has its roots in the construction by Gieseker of a degeneration of the moduli space of rank 2 vector bundles for a one-parameter family of smooth curves degenerating to a nodal curve \cite{gieseker-degeneration}, which was extended to higher rank by Nagaraj-Seshadri \cite{nagaraj-seshadri-I, nagaraj-seshadri-II}. Constructions of relative moduli spaces of semistable Gieseker vector bundles over $\overline{\cM}_{g,n}$ were obtained via Geometric Invariant Theory (GIT) by Caporaso \cite{caporaso_compactification_picard} in the rank 1 case, and Schmitt \cite{schmitt-compactification} for arbitrary rank. A Hitchin pair analog of the GIT arguments were developed by Balaji-Barik-Nagaraj 
\cite{balaji-barik-nagaraj} to obtain a relative moduli space of semistable Gieseker Hitchin pairs for a one-parameter family of smooth curves degenerating to a singular curve with a single node.

Our main strategy to obtain the universal compactifications in \Cref{introthm: moduli thm} over $\overline{\cM}_{g,n}$ is to bypass GIT and employ instead recent advances in the construction of moduli spaces via stacks by Alper-Halpern-Leistner-Heinloth \cite{alper-good-moduli, alper2019existence, halpernleistner2018structure}. We consider the moduli stack $\GHiggs_{g,n}^N$ of rank $N$ meromorphic Gieseker Higgs bundles, which admits a forgetful stabilization morphism
\[\GHiggs_{g,n}^N \to \overline{\cM}_{g,n}, \; \; \; (C, \sigma_{\bullet}, E, \psi) \mapsto (\overline{C}, \overline{\sigma}_{\bullet}).\]
In the unpunctured case (i.e. $n=0$), this stack of Gieseker Higgs bundles was recently considered in the work of Ben-Bassat-Das-Pantev \cite{benbassat-das-pantev}. Using the theory of $\Theta$-stability with respect to a numerical invariant as in the work of Heinloth and Halpern-Leistner \cite{heinloth-hilbertmumford, halpernleistner2018structure}, we are able to describe an open substack $\GHiggs_{g,n}^{N,ss} \subset \GHiggs_{g,n}^{N}$ of semistable meromorphic Gieseker Higgs bundles and show that it admits a relative good moduli space $\sqcup_{\chi \in \mathbb{Z}} \MGHiggs_{g,n}^{N,\chi}$ in the sense of Alper \cite{alper-good-moduli}, where each component $\MGHiggs_{g,n}^{N,\chi}$ is the space referred to in \Cref{introthm: moduli thm}. Our method of proof proceeds by using the strategy of infinite dimensional GIT developed in \cite{torsion-freepaper, gauged_theta_stratifications}, and the recent observation in the upcoming work \cite{gauged-maps-projective} that this strategy fits within the context of the moduli of Gieseker vector bundles.

Our other contribution in \Cref{introthm: moduli thm} is a study of the geometry of the stack $\GHiggs_{g,n}^N$ (see Subsection \ref{subsection: flatness of the hitchin morphism}). This allows us to prove directly that the corresponding Hitchin morphisms are syntomic at the level of stacks via a gluing argument.

It was observed by Das \cite{das-relative-log} that, in the unpunctured case ($n=0$) for a one-parameter family of curves as in the context of \cite{balaji-barik-nagaraj}, the symplectic structure on the generic fiber of the moduli of semistable Gieseker Higgs bundles naturally extends to a relative log-symplectic structure (in the sense of log-geometry). This was recently promoted, at the level of a derived enhancement of the stack $\GHiggs_{g,0} \to \overline{\cM}_{g,0}$, to a relative log $0$-shifted symplectic form by the work of Ben-Bassat-Das-Pantev \cite{benbassat-das-pantev}.

In the sequel to this paper \cite{dh-partII}, we will equip the stable loci in our moduli spaces $\MGHiggs_{g,n}^{N,\chi}$ and $\MGHiggs_{g,n}^{\cO_{\bullet},\chi}$ with naturally defined log-structures, and extend to
$\overline{\cM}_{g,n}$
the Poisson and symplectic structure over the boundary. We shall therefore complete the extension of the classical Poisson and symplectic systems attached to meromorphic Higgs bundles.

\noindent \textbf{Other future work and open problems.} Some of the motivation for our work comes from physics \cite{BDD-modified-hitchin-gln}, where at least the bases of the integrable systems are crucial for understanding superconformal field theories of class S. The physics is interested mostly in the nilpotent case; non-nilpotent 
deformations are ``massive", so not conformal. However, from a mathematical point of view, this restriction seems rather arbitrary. It would be interesting to carry out the extension to the case of non-nilpotent residues and to study the new phenomena that arise.

Another natural direction of research is the extension to meromorphic Higgs bundles for other reductive groups. Even though there is currently no good analogue of a notion of Gieseker $G$-bundle in the case when $G$ is an arbitrary reductive group, it is possible to attempt to study the geometry of the stack $\Higgs_{g,n}^{G}$ parametrizing meromorphic $G$-Higgs bundles on semistable twisted (stacky) curves. It would be interesting to carry out the analysis for $\Higgs_{g,n}^{G}$, construct the corresponding Hitchin morphisms, and define log-Poisson and log-symplectic forms at the level of the stacks.

\noindent \textbf{Acknowledgements.} We would like to thank Aswin Balasubramanian, Oren Ben-Bassat, Jacques Distler, Roberto Fringuelli, Daniel Halpern-Leistner, Johannes Horn and Tony Pantev for useful discussions related to this work.
During the preparation of this work, Ron Donagi was supported in part by NSF grants DMS 2001673 and 2401422, by NSF FRG grant DMS 2244978, and by Simons HMS Collaboration grant 390287. Both authors thank the Galileo Galilei Institute in Firenze, where some of this work was done.

\end{section}

\section{Meromorphic Gieseker Higgs bundles}
In this section, we introduce our main objects of interest: meromorphic Gieseker Higgs bundles. We define the stacks $\GHiggs_{g,n}^N$ that parametrize meromorphic Gieseker Higgs bundles on $n$-pointed curves, and prove some geometric properties of $\GHiggs_{g,n}^N$. We also construct a Hitchin morphism $H: \GHiggs_{g,n}^N \to \cA$ and show that it satisfies the existence part of the valuative criterion for properness.

\subsection{Some notation and preliminaries}
	We work over an algebraically closed field $k$ of characteristic $0$. We fix once and for all nonnegative integers $n,g$,  and denote by $\overline{\mathcal{M}}_{g,n}$ the Deligne-Mumford stack of stable $n$-pointed genus $g$ curves over $k$. 

 \begin{remark}
     Note that $\overline{\cM}_{g,n}$ is empty if $g=0$ and $n \leq 2$, or if $g=1$ and $n=0$. Our results are vacuous in those cases.
 \end{remark}

 Given a field extension $K \supset k$, recall that an $n$-pointed prestable curve $(C, \sigma_{\bullet})$ over $K$ is a geometrically connected projective curve $C$ over $K$ with at worst nodal singularities equipped with $n$-distinct points $\sigma_1, \sigma_2, \ldots, \sigma_n \in C(K)$ that lie in the smooth locus of $C$. Given a scheme $T$, a $T$-family of $n$-pointed prestable curves $(C, \sigma_{\bullet})$ is a flat proper finitely presented morphism $C \to T$ from an algebraic space $C$ and $n$-sections $\sigma_1, \sigma_2, \ldots, \sigma_n: T \to C$ such that for all $t \in T$, the fiber $C_t$ equipped with the points $\sigma_1(t), \sigma_2(t), \ldots, \sigma_n(t)$ is an $n$-pointed prestable curve over the residue field $k(t)$.

  \begin{notn}[Relative (log) canonical bundles]
     For any given scheme $T$ and a $T$-family of $n$-pointed prestable curves $(C, \sigma_{\bullet})$, we denote by $\omega_{C/T}$ the canonical line bundle on $C$ relative to $T$, and we set $\omega_{C/T}^{log} := \omega_{C/T}\left(\sum_{i=1}^n \sigma_i\right).$
 \end{notn}
 
 Recall that every family $(C, \sigma_{\bullet})$ of $n$-pointed prestable curves over a base scheme $T$ admits a stabilization $\varphi: (C, \sigma_{\bullet}) \to (\overline{C}, \overline{\sigma}_{\bullet})$, where $(\overline{C}, \overline{\sigma}_{\bullet})$ is a $T$-family of $n$-pointed stable curves. We say that $(C, \sigma_{\bullet})$ is a family of semistable $n$-pointed curves if $\varphi^*(\omega^{log}_{\overline{C}/T}) = \omega^{log}_{C/T}$. We denote by $\cM^{ss}_{g,n}$ the stack of semistable $n$-pointed genus $g$ curves, which is a smooth algebraic stack over $\Spec(k)$.  There is a stabilization morphism which is syntomic (i.e. a locally finitely presented flat morphism with lci fibers) of pure relative dimension $0$
	\[\stab: \cM^{ss}_{g,n} \to \overline{\mathcal{M}}_{g,n}, \; \; \; (C, \sigma_{\bullet}) \mapsto (\overline{C}, \overline{\sigma}_{\bullet}).\]

		\begin{defn}[Stack of coherent sheaves on stable curves]	We define $\Coh_{g,n}$ to be the pseudofunctor from $(\Sch/k)^{op}$ into groupoids that sends a test scheme $T$ to the groupoid of tuples $(C, \sigma_{\bullet}, \cF)$, where $(C, \sigma_{\bullet})$ is a $T$-family of $n$-pointed stable curves in $\overline{\cM}_{g,n}(T)$ and $\cF$ is a finitely presented $T$-flat $\cO_{C}$-module
  %modules 
  $\cF$.
		\end{defn}
	There is a forgetful morphism
 \[\Coh_{g,n} \to \overline{\cM}_{g,n}, \;\; \; (C, \sigma_{\bullet}, \mathcal{F}) \mapsto (C, \sigma_{\bullet})\]
 We note that we can equivalently think of $\Coh_{g,n}$ as the stack $\Coh_{\mathcal{C}/ \overline{\cM}_{g,n}}$ of coherent sheaves on the universal curve $\mathcal{C} \to \overline{\cM}_{g,n}$, and therefore $\Coh_{g,n}$ is an algebraic stack locally of finite type and with affine relative diagonal over $\overline{\mathcal{M}}_{g,n}$ \cite[\href{https://stacks.math.columbia.edu/tag/08WC}{Tag 08WC}]{stacks-project}. We will let $\Coh^{tf}_{g,n}$ denote the open substack of $\Coh_{g,n}$ parametrizing families $\mathcal{F}$ such that for all points $t \in T$ the restriction $\mathcal{F}|_{\mathcal{C}_{t}}$ is pure of dimension 1 (cf. \cite[Prop. 2.5]{torsion-freepaper}).

Let us record a useful fact, proven in \cite{gauged-maps-projective}.
\begin{prop} \label{prop: open closed generic rank}
    The substack $\Coh^{N}_{g,n} \subset \Coh^{tf}_{g,n}$ classifying families of torsion-free sheaves that have rank $N$ at every generic point of every fiber is an open and closed substack. \qed
\end{prop}

\subsection{The stack of Gieseker vector bundles}

In this subsection we explain some notation and results from the upcoming work \cite{gauged-maps-projective} which we use in this article.

\begin{defn}[Stack of Gieseker vector bundles, {\cite{gauged-maps-projective}}] \label{defn: gieseker vector bundles}
The stack $\GBun_{g,n}^{N}$ is the pseudofunctor from $(\Sch/k)^{op}$ into groupoids that sends a $k$-scheme $T$ to the groupoid of tuples $(C, \sigma_1, \ldots, \sigma_n, E)$ where
\begin{enumerate}[(1)]
	\item $(C, \sigma_{\bullet})$ is a family of $n$-pointed genus $g$ prestable curves over $T$ with stabilization morphism $\varphi: (C, \sigma_{\bullet}) \to (\overline{C}, \overline{\sigma}_{\bullet})$.
	\item $E$ is a rank $N$ vector bundle on $C$ such that $\det(E)$ is $\varphi$-ample, the counit $\varphi^* \varphi_*(E) \to E$ is surjective, and the pushforward $\varphi_*(E)$ is a $T$-flat family of torsion-free sheaves on $\overline{C} \to T$.
\end{enumerate}
\end{defn}

\begin{remark}
    It is not necessary in \Cref{defn: gieseker vector bundles} to impose the condition that $(C, \sigma_{\bullet})$ is semistable. Indeed, for any prestable $n$-pointed curve $(C, \sigma_{\bullet})$, the existence of a vector bundle $E$ as in (2) above forces the pointed curve to be semistable.
\end{remark}

\begin{remark} \label{remark: Gieseker bundle concrete}
It might help to have a concrete interpretation of what the conditions in \Cref {defn: gieseker vector bundles} say over a fixed stable curve $\overline{C}$ defined over the ground field $k$.
The curve $C$ is obtained from $\overline{C}$ by replacing each node of $ \overline{C}$ by a chain $R$ of $\ell \geq 0$ smooth rational curves $R_1,\dots,R_{\ell}$.
The restriction of the rank $N$ vector bundle $E$ to the chain $R$ is a direct sum $\bigoplus_{j=1}^N L_j$ of $N$ line bundles on $R$.
Every such line bundle $L$ on $R$ is obtained by gluing together line bundles $\cO(n_i)$ on $R_i \cong \mathbb{P}^1$ for some tuple of integers $n_{1}, \ldots n_{\ell}$.
Let us denote such a line bundle by $\cO(n_1, n_2, ..., n_{\ell})$, and set $L_j = \cO(n_{j,1}, n_{j,2}, \ldots, n_{j,\ell})$. 
The condition that $\varphi^* \varphi_*(E) \to E$ is surjective implies that all line bundles $L_j$ are globally generated on $R$, which means that each $n_{j,i}$ is nonnegative. 
The pushforward $\varphi_*(E)$ to the stable curve $\overline{C}$ would be locally free at the node $\varphi(R)$ if and only if all the $n_{j,i}$ vanished.
This pushforward is torsion-free at the node $\varphi(R)$ if and only if for each fixed value of $j$, the $n_{j,i}$ vanish, with at most one exception where $n_{j,i}=1$.
Relative ampleness of the determinant $\det(E)$ means that for each $i$, at least one $n_{j,i}$ is $>0$.

Altogether, Definition \ref {defn: gieseker vector bundles} says that the curve $C$ is obtained from  $\overline{C}$ by replacing each node of $\overline{C}$ by a chain of $\ell \leq N$ rational curves and taking $E$ to be the sum of line bundles $L_j = \cO(n_{j,1}, n_{j,2}, ..., n_{j,{\ell}})$, where each row of the $N \times \ell $ matrix $(n_{j,i})$ has all zeros except for at most a single $1$, while each column has at least one $1$. The pushforward $\varphi_*(E)$ is locally free if and only if all $(n_{j,i})=0$, which happens if and only if $\ell = 0$. This is the generic case. There is a finite number of possible non-zero matrices $(n_{j,i})$, and these describe ways for the vector bundle to degenerate to a torsion-free sheaf.
\end{remark}

  \begin{defn}[Cornalba line bundles] \label{defn: cornalba bundles}
		 For any given integer $m$, we define a line bundle $\mathcal{L}_{Cor,m}$ 
   on $\GBun_{g,n}^N$ as follows. For any $k$-scheme $T \to \overline{\cM}_{g,n}$ and any morphism
   $f: T \to \GBun_{g,n}^N$ corresponding to a tuple $(C, \sigma_{\bullet}, E)$, we set
			\[ f^*\mathcal{L}_{Cor,m} := \det\left(R(\pi_C)_*\left(\left(\omega_{C/T}^{log} \otimes \left(\det(E) \otimes \left(\omega_{\overline{C}/T}^{log}\right)^{\otimes m}|_{C}\right)^{\otimes 3} - [\mathcal{O}_C]\right)^{\otimes 2}\right)\right) \]
			where $\pi_C: C \to T$ is the structure morphism and $(\overline{C}, \overline{\sigma}_{\bullet}) \to T$ is the stabilization of $(C, \sigma_{\bullet})$. Here we take the Mumford-Knudsen determinant; note that the derived pushforward above is a perfect complex by  \cite[\href{https://stacks.math.columbia.edu/tag/0A1H}{Tag 0A1H}]{stacks-project}.
		\end{defn}

\begin{prop}[{\cite{gauged-maps-projective}}] \label{prop: gieseker vector bundles morphism proper}
    $\GBun_{g,n}^{N}$ is a smooth algebraic stack over $k$, equipped with a projective schematic morphism 
    \[\GBun_{g,n}^{N} \to \Coh^{N}_{g,n}, \; \; \; \; (C, \sigma_{\bullet}, E) \mapsto (\overline{C}, \overline{\sigma}_{\bullet}, \varphi_*(E)).\]
    The composition $\GBun_{g,n}^{N} \to \Coh^{N}_{g,n} \to \overline{\cM}_{g,n}$ is syntomic. \qed
\end{prop}

\begin{prop}[{\cite{gauged-maps-projective}}] \label{prop: gauged maps projective ampleness}
    Let $T$ be a quasicompact scheme equipped with a morphism $T \to \Coh^{N}_{g,n}$. Then for all sufficiently large $m\gg0$ the pullback of the line bundle $\mathcal{L}_{Cor,m}$ is $T$-ample on the proper scheme $\GBun_{g,n}^{N} \times_{\Coh^{N}_{g,n}} T$. \qed
\end{prop}

For each integer $\chi \in \mathbb{Z}$, there is an open and closed substack $\GBun_{g,n}^{N,\chi} \subset \GBun_{g,n}^{N}$ consisting of tuples $(C, \sigma_{\bullet}, E)$ such that the Euler characteristic of $E$ is $\chi$. This induces a decomposition $\GBun_{g,n}^{N} = \bigsqcup_{\chi \in \mathbb{Z}} \GBun_{g,n}^{N,\chi}$.

\subsection{Meromorphic Gieseker Higgs bundles}
\begin{defn}[Stack of torsion-free meromorphic Higgs sheaves]
	The stack $\Higgs^{tf}_{g,n}$ is the pseudofunctor from $(\Sch/k)^{op}$ to groupoids that sends a test scheme $T$ to the groupoid of tuples $(C, \sigma_{\bullet}, \cF, \psi)$, where $(C, \sigma_{\bullet}, \cF)$ is $T$-family in $\Coh^{tf}_{g,n}(T)$, and $\psi$ is a morphism of $\mathcal{O}_{C}$-modules $\psi: \mathcal{F} \to \mathcal{F} \otimes \omega_{C/T}^{log}$.
\end{defn}

\begin{notn}[Uniform rank $N$ meromorphic Higgs sheaves]
    We denote by $\Higgs^{N}_{g,n}$ the open and closed substack of $\Higgs^{tf}_{g,n}$ that parametrizes families $(C, \sigma_{\bullet}, \cF, \psi)$ such that $(C, \sigma_{\bullet}, \cF)$ belongs to $\Coh^N_{g,n}$.
\end{notn}

\begin{notn}
    For any given integer $\chi \in \mathbb{Z}$, we denote by $\Higgs^{N, \chi}_{g,n}$ the open and closed substack of $\Higgs^N_{g,n}$ that parametrizes torsion-free meromorphic Higgs sheaves whose underlying sheaf has Euler characteristic $\chi$.
\end{notn}

Recall that we can interpret $\Higgs_{g,n}^{tf}$ as the relative moduli stack $\Lambda\Coh^{tf}(\mathcal{C}/\overline{\cM}_{g,n})$ of torsion-free $\Lambda$-modules on the universal $n$-pointed curve $(\mathcal{C}, \sigma_{\bullet}) \to \overline{\cM}_{g,n}$, where $\Lambda := \Sym^{\bullet}_{\mathcal{O}_{\mathcal{C}}}\left((\omega_{\cC/\overline{\cM}_{g,n}}^{log})^{\vee}\right)$ (cf. \cite[Main example: Higgs bundles, pg. 86]{Simpson-repnI}). In particular, it is an algebraic stack locally of finite type and with affine relative diagonal over $\overline{\cM}_{g,n}$, by \cite[Prop. 2.23]{torsion-freepaper}. The same holds for the open and closed substack $\Higgs^{N}_{g,n}$.

\begin{defn}[Stack of rank $N$ meromorphic Gieseker Higgs bundles]
The stack $\GHiggs_{g,n}^{N}$ is the pseudofunctor from $(\Sch/k)^{op}$ into groupoids that sends a scheme $T$ to the groupoid of tuples $(C, \sigma_{\bullet}, E, \psi)$, where
\begin{enumerate}[(1)]
	\item $(C, \sigma_{\bullet}, E) \in \GBun_{g,n}^{N}(T)$.
	\item $\psi: E \to E\otimes \omega_{C/T}^{log}$ is a morphism of vector bundles over $C$.
\end{enumerate}
\end{defn}

\begin{lemma} \label{lemma: morphism from hitchin to bun is affine}
	The natural forgetful morphism 
	\[\GHiggs_{g,n}^{N} \to \GBun_{g,n}^{N}, \; \; \; (C, \sigma_{\bullet}, E, \psi) \mapsto (C, \sigma_{\bullet},E)\]
	is affine and of finite type. In particular, $\GHiggs_{g,n}^{N}$ is an algebraic stack locally of finite type over $\overline{\cM}_{g,n}$
 and with affine relative diagonal over $\overline{\cM}_{g,n}$. 
\end{lemma}
\begin{proof}
	Let $T$ be a scheme, and let $T \to \GBun_{g,n}^{N}$ be a morphism corresponding to a tuple $(C, \sigma_{\bullet}, E)$. The fiber product $\GHiggs_{g,n}^{N} \times_{\GBun_{g,n}^{N}} T$ is the functor that sends a scheme $Y \to T$ to the set of global sections $H^0\left(\left(\mathcal{E}nd(E) \otimes \omega_{C/T}^{log}\right)|_{C_{Y}}\right)$. Since the sheaf $\mathcal{E}nd(E)\otimes \omega_{C/T}^{log}$ on the $T$-proper scheme $C \to T$ is $T$-flat and finitely presented, this functor is represented by an affine scheme of finite presentation over $T$ \cite[\href{https://stacks.math.columbia.edu/tag/08K6}{Tag 08K6}]{stacks-project}.
\end{proof}

\begin{notn}
    For any given integer $\chi \in \mathbb{Z}$, we denote by $\GHiggs_{g,n}^{N,\chi}$ the open and closed substack of $\GHiggs_{g,n}^N$ given by the preimage of $\GBun_{g,n}^{N,
    \chi} \subset \GBun_{g,n}^N$ under the morphism $\GHiggs_{g,n}^N \to \GBun_{g,n}^N$.
\end{notn}

\begin{lemma} \label{lemma: properness of gieseker hitchin stack over torsion-free hitchin stack}
	The pushforward to the stabilization induces a proper schematic morphism of stacks
	\[\varpi: \GHiggs_{g,n}^{N} \to \Higgs^{N}_{g,n}, \; \; (C, \sigma_{\bullet}, E, \psi) \mapsto (\overline{C}, \overline{\sigma}_{\bullet}, \varphi_*(E), \varphi_*(\psi)).\]
\end{lemma}
\begin{proof}
	Consider the fiber product $F:=\GBun_{g,n}^{N} \times_{\Coh^{N}_{g,N}}\Higgs^{N}_{g,n}$, which fits into the following commutative diagram
\[\begin{tikzcd}
				\GHiggs_{g,n}^{N} \ar[rr, "i"] \ar[dr, "\varpi"] &  & F \ar[dl]\\   & \Higgs^{N}_{g,n} &
			\end{tikzcd}\]
	By \Cref{prop: gieseker vector bundles morphism proper}, the morphism $F=\GBun_{g,n}^{N} \times_{\Coh^{tf}_{g,n}} \Higgs^{N}_{g,n} \to \Higgs^{N}_{g,n}$ is schematic and proper. We shall prove that the horizontal morphism $i$ is a closed immersion, thus concluding the proof.
	Let $T$ be a scheme, and let $T \to F$ be a morphism corresponding to a tuple $(C, \sigma_{\bullet}, E, \psi)$, where $(C, \sigma_{\bullet}, E)$ is a tuple in $\GBun_{g,n}^{N}(T)$ with stabilization morphism $\varphi: (C,\sigma_{\bullet}) \to (\overline{C}, \overline{\sigma}_{\bullet})$, and $\psi$ is a morphism $\psi: \varphi_*(E) \to \varphi_*(E) \otimes \omega_{\overline{C}/T}^{log}$ of quasicoherent sheaves on $\overline{C}$. Note that, since the $T$-family of $n$-pointed curves $(C, \sigma_{\bullet})$ is semistable, there is a natural identification $\varphi^*(\omega_{\overline{C}/T}^{log}) \cong \omega_{C/T}^{log}$. The fiber product $\GHiggs_{g,n}^{N} \times_{F} T$ is the functor that sends a scheme $Y \to T$ into the set of morphisms $\theta: E|_{C_{Y}} \to E|_{C_{Y}}\otimes\omega_{C_Y/Y}^{log}$ such that the pushforward $(\varphi_{Y})_*(\theta)$ under the base-changed stabilization morphism $\varphi_{Y}: C_{Y} \to \overline{C}_{Y}$ corresponds to the base-change $\psi_Y:= \psi|_{\overline{C}_Y}: \varphi_*(E)|_{\overline{C}_{Y}} \to \varphi_*(E)|_{\overline{C}_{Y}} \otimes \omega_{\overline{C}_Y/Y}^{log}$ under the natural identifications. Such $\theta$ must fit into the following commutative diagram, where the vertical morphisms are the surjective counits
\[\begin{tikzcd}
		(\varphi_{Y})^*(\varphi_{Y})_*(E) \ar[r, "(\varphi_Y)^*(\psi_{Y})"] \ar[d] &(\varphi_{Y})^*(\varphi_Y)_*(E) \otimes \omega_{C_Y/Y}^{log} \ar[d] \\ E|_{C_{Y}} \ar[r, "\theta"] &  E|_{C_{Y}} \otimes \omega_{C_Y/Y}^{log}.
	\end{tikzcd}\]
In particular such a morphism exists if and only if the composition 
\[(\varphi_{Y})^*(\varphi_{Y})_*(E) \xrightarrow{\varphi_{Y}^*(\psi_{Y})} (\varphi_{Y})^*(\varphi_Y)_*(E) \to E|_{C_{Y}} \otimes \omega_{C_Y}^{log} \]
factors through $E|_{C_{Y}}$, and if that is the case then $\theta$ is completely determined. If we denote by $\mathcal{K}_{Y}$ the kernel of the surjective counit $(\varphi_{Y})^*(\varphi_{Y})*(E) \to E|_{C_{Y}}$, then it follows that $\GHiggs_{g,n}^{N} \times_{F} T$ is the subfunctor of $T$ consisting of those morphisms $Y \to T$ such that the composition 
\[g_{Y}: \mathcal{K}_{Y} \to (\varphi_{Y})^*(\varphi_{Y})_*(E) \xrightarrow{\varphi_{Y}^*(\psi_{Y})} (\varphi_{Y})^*(\varphi_Y)_*(E) \to E|_{C_{Y}} \otimes \omega_{C_Y/Y}^{log} \]
vanishes. But, since $E$ is flat over $T$, it follows that $\mathcal{K}_{Y} = \mathcal{K}_{T}|_{C_{Y}}$, and the construction implies that $g_{Y}: \mathcal{K}_{Y} \to E|_{C_{Y}} \otimes \omega_{C_Y/Y}^{log}$ is the pullback of the morphism $g_{T}: \mathcal{K}_{T} \to E \otimes \omega_{C/T}^{log}$ under $C_{Y} \to C$. In other words, if we denote by $\Hom_{C/T}\left(\mathcal{K}_T, \, E \otimes \omega_{C/T}^{log}\right)$ the functor that sends a scheme $Y \to T$ into the set of morphisms $\Hom_{C_{Y}}\left(\mathcal{K}_{Y}, \, E|_{C_{Y}} \otimes \omega_{C/T}^{log}|_{C_{Y}}\right)$, then we have the following Cartesian diagram
\[\begin{tikzcd}
		\GHiggs_{g,n}^{N} \times_{F} T \ar[r] \ar[d, "\text{pr}_2"] &  \ar[d, "0"] T \\ T \ar[r, "g_T"] &  \Hom_{C/Y}\left(\mathcal{K}_{T}, E \otimes\omega_{C/T}^{log}\right)
	\end{tikzcd}\]
By \cite[\href{https://stacks.math.columbia.edu/tag/08K6}{Tag 08K6}]{stacks-project}, the natural morphism $\Hom_{C/T}\left(\mathcal{K}_T, \, E \otimes \omega_{C/T}^{log}\right) \to T$ is affine, and so in particular the $0$-section $0: T \to \Hom_{C/T}\left(\mathcal{K}_T, \, E \otimes \omega_{C/T}^{log}\right)$ is a closed immersion. It follows that the base-change $\GHiggs_{g,n}^{N} \times_{F} T \xrightarrow{\text{pr}_2} T$ is a closed immersion.
\end{proof}

\subsection{The Hitchin morphism}
 \begin{defn}[Hitchin base]
 	The Hitchin base $\cA$ is the functor from $(\Sch/k)^{op}$ into groupoids that sends a scheme $T$ into the groupoid of families of $n$-pointed stable curves $(C ,\sigma_{\bullet})$ in $\overline{\cM}_{g,n}(T)$ along with tuples $(a_i)_{i=1}^N$, where $a_i \in H^0\left(C, (\omega_{C/T}^{log})^{\otimes i}\right)$. 
 \end{defn}

 \begin{lemma} \label{lemma: meromorphic hitchin base vector bundle}
	The forgetful morphism 
 \[\cA \to \overline{\cM}_{g,n}, \; \; \; (C, \sigma_{\bullet}, a_{\bullet}) \mapsto (C, \sigma_{\bullet})\]
 is affine and of finite type. Furthermore, $\cA$ is naturally isomorphic to the total space of a vector bundle on $\overline{\cM}_{g,n}$.
\end{lemma}
\begin{proof}
\cite[\href{https://stacks.math.columbia.edu/tag/08K6}{Tag 08K6}]{stacks-project} implies that $\cA \to \overline{\cM}_{g,n}$ is affine and of finite type. In order to prove that $\cA$ is a vector bundle, by cohomology and base-change and the fact that $\overline{\cM}_{g,n}$ is reduced, it suffices to show that the dimension of $H^0(C, (\omega_{C}^{log})^{{\otimes}i})$ is constant for all $n$-pointed genus $g$ stable curves defined over a field. By Riemann-Roch, this dimension is $(2g-2+n)\cdot i + 1-g$ when $n>0$ or $i>1$, and it is equal to $g$ when $n=0$ and $i=1$. Regardless of the value of $i$ and $n$, we see that the dimension is independent of the chosen $n$-pointed stable curve.
\end{proof}

Let $T$ be a scheme and let $T \to \GHiggs_{g,n}^{N}$ be a morphism represented by a tuple $(C, \sigma_{\bullet}, E, \psi)$, with stabilization $\varphi: (C, \sigma_{\bullet}) \to (\overline{C}, \overline{\sigma}_{\bullet})$. For each $1 \leq i \leq N$ we have an induced wedge morphism 
\[\wedge^i \psi : \bigwedge^i E \to \bigwedge^i(E \otimes \omega_{C/T}^{log}) \cong \bigwedge^i(E) \otimes (\omega_{C/T}^{log})^{\otimes i}.\]
Using the trace homomorphism $\tr: \mathcal{E}nd(\bigwedge^iE) \to \mathcal{O}_{C}$, we get an element 
\[a_i(\psi):=\tr\otimes\text{id}_{(\omega_{C/T}^{log})^{\otimes i}}(\wedge^i\psi) \in H^0\left(C, (\omega_{C/T}^{log})^{\otimes i}\right).\]
We can think of $a_i(\psi)$ alternatively as an element 
$\overline{a}_i \in H^0\left(\overline{C}, \varphi_{*}((\omega_{C/T}^{log})^{\otimes i})\right)$. Since $(C, \sigma_{\bullet})$ is a family of semistable curves, there is a natural identification $\varphi^*(\omega_{\overline{C}/T}^{log}) \cong \omega_{C/T}^{log}$. Using the projection formula and the fact that $\varphi_*(\cO_{C}) = \cO_{\overline{C}}$, we conclude that there is a natural identification $\varphi_{*}((\omega_{C/T}^{log})^{\otimes i}) \cong (\omega_{\overline{C}/T}^{log})^{\otimes i}$, and therefore we may view $\overline{a}_i(\psi)$ as an element of $ H^0\left(\overline{C}, (\omega_{\overline{C}/T}^{log})^{\otimes i}\right)$.

This yields the following morphism of $\overline{\cM}_{g,n}$-stacks,
\[H: \GHiggs_{g,n}^{N} \to \cA, \; \; \; (C, \sigma_{\bullet}, E,\psi) \mapsto (\overline{C}, \overline{\sigma}_{\bullet}, (\overline{a}_i)_{i=1}^N).\]
\begin{notn}[Hitchin morphism]
    We call $H: \GHiggs_{g,n}^{N} \to \cA$ the (universal) Hitchin morphism.
\end{notn}

Let $(\cC, \sigma_{\bullet}, a_{\bullet}) \to \cA$ denote the universal family over $\cA$. Consider the total space $p: \Tot(\omega_{\cC/\cA}^{log}) \to \cC$ of the line bundle $\omega_{\cC/\cA}^{log}$. There is a universal global section $\lambda \in H^0\left(\Tot(\omega_{\cC/\cA}^{log}), p^*(\omega_{\cC/\cA}^{log})\right)$, which we can use to define a section
\[s = \lambda^{\otimes N} + \sum_{i=1}^{N}(-1)^i p^*(a_i) \otimes \lambda^{\otimes (N-i)}\]
in $H^0\left(\Tot(\omega_{\cC/\cA}^{log}), p^*(\omega_{\cC/\cA}^{log})^{\otimes N}\right)$. The vanishing locus $\widetilde{\cS} \hookrightarrow \Tot(\omega_{\cC/\cA}^{log})$ of the section $s$ is a relative effective Cartier divisor on $p: \Tot(\omega_{\cC/\cA}^{log}) \to \cC_{\cA}$, and the induced map $\widetilde{\cS} \to \cC$ is a finite flat morphism of degree $N$. 

\begin{notn}
    We call $\widetilde{\cS} \to \cC$ the universal spectral curve. For each morphism $T \to \cA$ we call the base change $\widetilde{\cS}_{T} \to \cC_{T}$ the spectral curve associated to the spectral data $T \to \cA$.
\end{notn}

Let $\Coh^{tf}_{\widetilde{\cS}/\cA}$ denote the relative stack of torsion-free sheaves on the family of projective curves $\widetilde{\cS} \to \cA$. As in the usual BNR correspondence \cite{bnr, simpson-repnII}, pushing forward under the finite flat morphism $\widetilde{\cS} \to \cC$ induces a monomorphism of stacks $\Coh^{tf}_{\widetilde{\cS}/\cA} \hookrightarrow \Higgs^{tf}_{g,n}\times_{\overline{\cM}_{g,n}} \cA$.
\begin{prop}[BNR for Gieseker Higgs bundles] \label{prop: bnr} The morphism 
\[\varpi \times_{\overline{\cM}_{g,n}} H: \GHiggs_{g,n}^{N} \to \Higgs^{tf}_{g,n} \times_{\overline{\cM}_{g,n}} \cA\]
(where $\varpi$ is as in \Cref{lemma: properness of gieseker hitchin stack over torsion-free hitchin stack})
factors through $\Coh^{tf}_{\widetilde{\cS}/\cA} \hookrightarrow \Higgs^{tf}_{g,n} \times_{\overline{\cM}_{g,n}} \cA$, thus inducing a uniquely determined morphism $c:\GHiggs_{g,n}^{N} \to \Coh^{tf}_{\widetilde{\cS}/\cA}$.
\end{prop}
\begin{proof}
    Let $T$ be a scheme equipped with a morphism $T \to \cA$. Let $Y$ be a $T$-scheme, and let $(C, \sigma_{\bullet}, E, \psi)$ denote a $Y$-point of $\GHiggs_{g,n}^{N} \times_{\cA} T$. Let $(\cF, \tau)$ denote the corresponding $Y$-point of $\Higgs^{tf}_{g,n} \times_{\overline{\cM}_{g,n}} T$ obtained by pushing forward to the stabilization $\overline{C}_Y$. The pair $(\cF, \tau)$ is equivalent to a $Y$-family $\cG$ of coherent pure sheaves with relative proper support on the total space $\Tot(\omega_{\overline{C}_Y/Y}^{log}) \to \overline{C}_Y \to Y$. It suffices to show that $\cG$ is supported on the corresponding spectral curve $\widetilde{\cS}_Y \subset \Tot(\omega_{\overline{C}_Y/Y}^{log})$ determined by the spectral data $Y \to T \to \cA$. To simplify notation, we assume without loss of generality that $T=Y$. 
    
    By the usual BNR correspondence applied to the pair $(E, \psi)$ on $C$, there is a $T$-family of torsion-free sheaves $\cG^C$ on the spectral curve $\widetilde{\cS}_T \times_{\overline{C}} C \to C \to T$ such that $(E, \psi)$ is obtained from the pushforward of $\cG^C$ under $\widetilde{\cS}_T \times_{\overline{C}} C \to C$. Moreover, it follows by considering the pushforward under the commutative square
\[
	\begin{tikzcd}
		\widetilde{\cS}_T \times_{\overline{C}} C \ar[r] \ar[d] &  \ar[d] \widetilde{\cS}_T \\ C \ar[r] &  \overline{C}
	\end{tikzcd}
\]
    that $\cG$ is the pushforward of $\cG^C$ under the morphism $\widetilde{\cS}_T \times_{\overline{C}} C \to \widetilde{\cS}_T \hookrightarrow \Tot(\omega_{\overline{C}_T/T}^{log}))$. Therefore, $\cG$ is supported on $\widetilde{\cS}_T$, as desired.
\end{proof}
The following was shown in \cite[\S 5]{benbassat-das-pantev} in the case when $n=0$. We provide an alternative proof by appealing to the results in \cite{gauged-maps-projective}.
\begin{prop} \label{prop: properness of the hitchin morphism}
    The Hitchin morphism $H: \GHiggs_{g,n}^{N} \to \cA$ satisfies the existence part of the valuative criterion of properness for stacks.
\end{prop}
\begin{proof}
        Fix a discrete valuation ring $R$ and a morphism $\Spec(R) \to \mathcal{A}$. Let $\eta$ denote the generic point of $\Spec(R)$, and fix an $\eta$-point $(C, \sigma_{\bullet}, E, \psi)$ of $\GHiggs_{g,n}^{N}(\eta)$ lying over $\eta \to \Spec(R) \to \cA$. We denote by $(\overline{C}, \overline{\sigma}_{\bullet}) \to \Spec(R)$ the $n$-pointed stable curve corresponding to  $\Spec(R) \to \mathcal{A} \to \overline{\cM}_{g,n}$, and $\widetilde{C} \to \overline{C}$ the corresponding spectral curve. Let $(\cF, \tau)$ denote the torsion-free meromorphic Higgs bundle on $\overline{C}_{\eta}$ obtained by pushing forward $(E, \psi)$ to the stabilization. By \Cref{prop: bnr}, the meromorphic Higgs sheaf $(\cF, \tau)$ comes from the pushforward of a torsion-free sheaf $\cG$ on $\widetilde{C}_{\eta}$. Since the morphism $\Coh^{tf}_{\widetilde{C}/\Spec(R)} \to \Spec(R)$ satisfies the existence part of the valuative criterion for properness \cite[Prop. 6.16]{torsion-freepaper}, we can extend $\cG$ to a torsion-free sheaf on $\widetilde{C}$, and by pushing forward we obtain a torsion-free meromorphic Higgs sheaf on $\overline{C}$ extending $(\cF, \tau)$. In other words, we can extend $(\cF, \tau)$ to an $R$-point of $\Higgs^{N}_{g,n} \times_{\overline{\cM}_{g,n}} \Spec(R)$. Since the pushforward morphism $\GHiggs_{g,n}^{N} \to \Higgs^{N}_{g,n}$ is schematic and proper, we can further lift this to an $R$ point of $\GHiggs_{g,n}^{N} \times_{\overline{\cM}_{g,n}} \Spec(R)$ extending $(\overline{C}, \overline{\sigma}_{\bullet})$. Since $\cA \to \overline{\cM}_{g,n}$ is separated, this extension actually factors through $\GHiggs_{g,n}^{N} \times_{\cA} \Spec(R)$, as desired.
\end{proof}

\subsection{General fibers of the Hitchin morphism}

In this subsection we give a description of the fibers of the Hitchin morphism over the following open locus.
\begin{defn}[A locus of nodal spectral curves]
We denote by $\cU \subset \cA$ the open substack consisting of spectral data $T \to \cA$ such that the corresponding spectral cover $\widetilde{\cS}_T \to \cC_T$ satisfies all the following
\begin{enumerate}[(1)]
    \item The $T$-fibers of $\widetilde{\cS}_T$ are nodal.
    \item For each $t \in T$, the morphism $\widetilde{\cS}_t \to \cC_t$ is \'etale at every point that maps to a node of $\cC_t$.
    \item For each $t \in T$, the smooth locus of $\widetilde{\cS}_t$ is the preimage of the smooth locus of $\cC_t$.
\end{enumerate}
\end{defn}

Over the open locus $\cU \subset \cA$, the fibers of the morphism $H: \GHiggs_{g,n}^{N} \to \cA$ are heavily related to the stack of torsion-free sheaves on the corresponding spectral curve.
\begin{prop}
    The restriction $c: \GHiggs_{g,n}^{N} \times_{\cA} \cU \to \Coh^{tf}_{\widetilde{\cS}_{\cU}/\cU}$ of the morphism defined in \Cref{prop: bnr} over the open $\cU \subset \cA$
    is a finite morphism. The image of this morphism is the open and closed substack $\Coh^{1}_{\widetilde{\cS}_{\cU}/\cU} \subset \Coh^{tf}_{\widetilde{\cS}_{\cU}/\cU}$
    consisting of torsion-free sheaves on the (nodal curve) fibers of the family $\widetilde{\cS}_{\cU} \to \cU$ that have rank 1 at all generic points of irreducible components (cf. \Cref{prop: open closed generic rank}).
\end{prop}
\begin{proof}
    Since the composition $\GHiggs_{g,n}^{N}\times_{\cA} \cU \to \Coh^{tf}_{\widetilde{\cS}_{\cU}/\cU} \hookrightarrow \Higgs^{tf}_{g,n}\times_{\overline{\cM}_{g,n}} \cU$ is schematic and proper, it follows that $c: \GHiggs_{g,n}^{N} \times_{\cA} \cU \to \Coh^{tf}_{\widetilde{\cS}_{\cU}/\cU}$ is a proper schematic morphism. By construction, the morphism $c$ factors through the open and closed substack $\Coh^{ 1}_{\widetilde{\cS}_{\cU}/\cU}$ as in the statement of the proposition. It suffices to show that the induced morphism $c: \GHiggs_{g,n}\times_{\cA} \cU \to \Coh^{ 1}_{\widetilde{\cS}_{\cU}/\cU}$ is quasifinite, since we already know that it is proper and schematic.

Fix an algebraically closed field $K \supset k$ and a morphism $\Spec(K) \to \cU$. This yields the data of an $n$-pointed stable curve $(\overline{C}, \overline{\sigma}_{\bullet})$ over $K$ and a corresponding nodal spectral cover $\widetilde{\cS}_K \to \overline{C}$ of degree $N$. Consider the stack $\cM^{ss}_{\overline{C}}:= \cM^{ss}_{g,n} \times_{\overline{\cM}_{g,n}} \Spec(K)$ of all $n$-pointed semistable curves with stabilization $(\overline{C}, \overline{\sigma}_{\bullet})$. The universal curve $\cC^{ss} \to \cM_{\overline{C}}^{ss}$ admits a (pointed) stabilization morphism $\cC^{ss} \to \overline{C} \times_{K} \cM_{\overline{C}}^{ss}$. We can take the fiber product with the spectral cover $\widetilde{\cS}_K \times_{K} \cM_{\overline{C}}^{ss}  \to \overline{C} \times_{K} \cM_{\overline{C}}^{ss}$ to obtain a universal base-changed spectral cover $\widetilde{\cS^{ss}} \to \cC^{ss}$ over the universal semistable modification of $\overline{C}$. The usual BNR correspondence identifies $\GHiggs_{g,n} \times_{\cA} \Spec(K)$ with an open substack of the relative stack $\Coh^{1}_{\widetilde{\cS^{ss}}/ \cM^{ss}_{\overline{C}}}$ of torsion-free sheaves of uniform rank 1 on the fibers of the family of nodal curves $\widetilde{\cS^{ss}} \to \cM^{ss}_{\overline{C}}$. Namely, it consists of the open locus of sheaves whose pushforward under $\widetilde{\cS^{ss}} \to \cC^{ss}$ are Gieseker vector bundles. By the definition of $\cU$, the morphism $\widetilde{\cS^{ss}} \to \cC^{ss}$ is \'etale over the nodes of each $\cM^{ss}_{\overline{C}}$-fiber of $\cC^{ss}$. By working \'etale locally on fibers, then it can be seen that if the pushforward of an element in $\Coh^{1}_{\widetilde{\cS^{ss}}/ \cM^{ss}_{\overline{C}}}$ is a vector bundle of rank $N$ on a fiber of $\cC^{ss}$, then this forces the sheaf on the spectral cover (fiber of $\widetilde{\cS^{ss}}$) to be a line bundle. 

Therefore, we may view $\GHiggs_{g,n}\times_{\cA} \Spec(K)$ as an open substack of the relative Picard stack $\cP ic_{\widetilde{\cS^{ss}}/ \cM^{ss}_{\overline{C}}}$, where we only consider certain ``allowed" semistable modifications $\Sigma \to \widetilde{\cS}_K$ of the spectral cover $\widetilde{\cS}_K$. An ``allowed" semistable modification of $\widetilde{\cS}_K$ looks as follows. For each node of $\overline{C}$, the preimage in $\widetilde{\cS}_K$ consists of $N$ nodes. This way, we can partition the set of nodes of $\widetilde{\cS}_K$ into equivalence classes of $N$-nodes. A semistable modification of $C \to \overline{C}$ consists of adding chains of rational curves separating some of the nodes in $\overline{C}$. The corresponding modification $\Sigma \to \widetilde{\cS}_K$ given by the fiber product
  \[
	\begin{tikzcd}
		\Sigma \ar[r] \ar[d] & \widetilde{\cS}_K \ar[d] \\ C \ar[r] &   \overline{C}
	\end{tikzcd}
\]
adds chains of rational curves on the corresponding preimages in $C$. So ``allowed" semistable modifications are just semistable modifications where the corresponding chain of rational curves have the same length for any two nodes in the same equivalence class.

For any such modification $\Sigma \to \widetilde{\cS}_K$, we have the cover morphism $\Sigma \to C$. A line bundle $\cL$ on $\Sigma$ represents an element of the stack $\GHiggs_{g,n}\times_{\cA} \Spec(K)$ if and only if its pushforward to $C$ is a Gieseker vector bundle. This condition can be described concretely for $\cL$ as follows. For each connected chain $R \subset C$ of rational curves in $C$ contracted by $C \to \overline{C}$, consider the $N$ distinct connected chains of rational curves $\{R_i\}_{i=1}^N$ in $\Sigma$ that map to $R$ under $\Sigma \to C$. Suppose that the chain $R$ has length $l$, and so it consists of $l$ copies of $\mathbb{P}^1$ that we can order and label by indexes $1 \leq j \leq l$ such that the $j^{th}$ copy of $\mathbb{P}^1$ is attached to the $(j-1)^{th}$ and the $(j+1)^{th}$ copies. For each $j$ between $1$ and $l$, we denote by $\cH_j$ the unique line bundle on $R$ that restricts to $\cO_{\mathbb{P}^1}(1)$ on the $j^{th}$ copy of $\mathbb{P}^1$ and restricts the trivial line bundle on all other copies. By convention, let us set $\cH_0$ to be the trivial line bundle on $R$. Since the morphism $\Sigma \to C$ induces a canonical isomorphism of each $R_i$ with $R$, we can define line bundles $\cH_j$ on each $R_i$ via pullback. For each $N$-tuple of numbers $\vec{i} = (i_1, i_2, \ldots, i_N)$ satisfying $0 \leq i_j \leq l$, we denote by $\cH_{\vec{i}}$ the unique line bundle on the union $\sqcup R_i$ that restricts to $\cH_{i_j}$ on the chain $R_{i}$. Then the line bundle $\cL$ on $\Sigma$ represents a point in $\GHiggs_{g,n}\times_{\cA} \Spec(K)$ if and only if its restriction $\cL|_{\sqcup_i R_i}$ is isomorphic to $\cH_{\vec{i}}$ for some tuple $\vec{i}$ such that all numbers between $1$ and $l$ appear at least ones among the tuple $i_j$. Roughly, this says that in each chain $R_i$, the restriction is either trivial or has a single rational curve in the chain where it becomes $\cO_{\mathbb{P}^1}(1)$, and that for each index $1 \leq j \leq l$ there is at least one chain $R_i$ where the line bundle restricts to $\cO_{\mathbb{P}^1}(1)$ on the $\mathbb{P}^1$ at index $j$.

  Note that, by contracting the rational curves on the chains where the line bundle $\cL$ is trivial, we obtain partial contraction $\Sigma \to \overline{\Sigma} \to \widetilde{\cS}_K$ such that $\overline{\Sigma} \to \widetilde{\cS}_K$ is a quasistable modification as in Caporaso's work \cite[\S3.3, Defn., pg. 616]{caporaso_compactification_picard}, and the pushforward of $\cL$ to $\overline{\Sigma}$ is a line bundle that restricts to $\cO_{\mathbb{P}^1}(1)$ on each rational curve contracted by $\overline{\Sigma} \to \widetilde{\cS}_K$.   
     
     Using this description, we prove now quasifiniteness of $c_K$. Given a uniform rank 1 torsion-free sheaf $\cF$ on $\widetilde{\cS}_K$, there is a uniquely determined isomorphism class of quasistable modification $\overline{\Sigma} \to \widetilde{\cS}_K$ along with a line bundle $\cK$ on $\overline{\Sigma}$ such that it restricts to $\cO_{\mathbb{P}^1}(1)$ on each contracted component of $\overline{\Sigma} \to \widetilde{\cS}_K$ and such that its pushforward to $\widetilde{\cS}_K$ recovers $\cF$ \cite[Thm. 6.1]{esteves-pacini}. By our discussion above, the possible preimages of $\cF$ in $\GHiggs_{g,n}\times_{\cA} \Spec(K)$ consists of the set of choices of how to add further rational curves to obtain an ``allowed" semistable modification
     \[
	\begin{tikzcd}
		\Sigma \ar[r] \ar[d] & \overline{\Sigma} \ar[r]& \widetilde{\cS}_K \ar[d] \\ C \ar[rr] & &  \overline{C}
	\end{tikzcd}
\]
such that the pullback 
$\cK|_{\Sigma}$ satisfies the conditions that we described when restricted to each set of chains 
$\sqcup_i R_i$ (note that the curve 
$\Sigma$ determines completely the isomorphism class of $C$, so we can ignore $C$ as part of the data). There are only finitely many choices of how to add rational curves to obtain $\Sigma$ (note that the possible lengths of chains of rational curves is always bounded by $N$), and so there are only finitely many preimages of $\cF$ in $\GHiggs_{g,n}\times_{\cA} \Spec(K)$.
\end{proof}

% We observe by means of an example that the comparison morphism $c: \GHiggs_{g,n}^{N} \to \Coh^{tf}_{\widetilde{\cS}_{\cA}/\cA}$ can have fibers of positive dimension over the points in $\cA$ that correspond to more ``degenerate" spectral curves.
% \begin{example}
%     \andres{TO DO}
% \end{example}

\begin{section}{Residues}

We may think of the Hitchin morphism $H: \GHiggs^N_{g,n} \to \cA$ as a log-Poisson system relative to $\overline{\cM}_{g,n}$, which is furthermore log-symplectic if $n=0$. This is in the sense of log-geometry, as will be made precise in \cite{dh-partII}; see \cite{benbassat-das-pantev} for the case when $n=0$.
The purpose of this section is to describe the log-symplectic leaves of this system. Similarly to the classical case of a fixed smooth projective curve, they will consist of meromorphic Gieseker Higgs bundles with a fixed given set of residue classes at the markings. 

We focus on the case when all the residues are nilpotent. For any given $n$-tuple of nilpotent conjugacy classes $\cO_{\bullet}$, we define a stack $\GHiggs_{g,n}^{\cO_{\bullet}}$ of meromorphic Gieseker Higgs bundles with residue classes given by the tuple $\cO_{\bullet}$. We show that $\GHiggs_{g,n}^{\cO_{\bullet}}$ is an equidimensional local-complete-intersection stack, and we compute its dimension. Finally, we construct a syntomic Hitchin morphism $H: \GHiggs_{g,n}^{\cO_{\bullet}} \to \cA^{\cO_{\bullet}}$ to the modified Hitchin base introduced by Balasubramanian-Distler-Donagi \cite{BDD-modified-hitchin-gln}.

\subsection{The residue morphism}
Let $T$ be a scheme, and let $(C, \sigma_{\bullet}) \to T$ be a family of $n$-pointed semistable curves in $\cM^{ss}_{g,n}(T)$. Since each marking $\sigma_i: T \to C$ lies in the smooth locus of $C \to T$, there is a canonical relative residue isomorphism $r_i: \sigma_i^*\omega_{C/T}^{log} \xrightarrow{\sim} \cO_T$.

In the following definition, we consider the quotient stack $\mathfrak{gl}_N/\GL_N$ where $\GL_N$ acts via its adjoint representation. A $T$-point of $\mathfrak{gl}_N/\GL_N$ consists of a pair $(E, \psi)$ where $E$ is a rank $N$ vector bundle on $T$ and $\psi$ is an endomorphism $\psi: E \to E$.

\begin{remark}
    We may think of $(E,\psi)$ as a family of Higgs bundles on (a family of) points over $T$, instead of a family of curves over $T$.
\end{remark}

\begin{defn}[Residue morphism]
    We denote by $res_i: \GHiggs_{g,n}^N \to \mathfrak{gl}_N / \GL_N$ the morphism of stacks that for any $k$-scheme $T$ sends a tuple $(C, \sigma_{\bullet}, E, \psi)$ to the pair consisting of the vector bundle $\sigma_i^*(E)$ and the endomorphism 
    \[r_i(\psi): \sigma_i^*(E) \xrightarrow{\sigma^*_i(\psi)} \sigma_i^*(E) \otimes \sigma_i^*\omega_{C/T}^{log} \xrightarrow{id \otimes r_i} \sigma_i^*(E).\]
    We denote by $\vec{res}$ the product $\prod_i res_i: \GHiggs_{g,n}^N \to \left(\mathfrak{gl}_N/\GL_N\right)^n$.
\end{defn}

The Fuchs relations force the morphism $\vec{res}$ to factor through a closed substack $F_N \subset \left( \mathfrak{gl}_N/\GL_N\right)^N$ which we define next.
\begin{defn}[Fuchs stack]
    The Fuchs stack $\mathscr{F}_N$ is the closed substack of $\left( \mathfrak{gl}_N/\GL_N\right)^n$ that sends a $k$-scheme $T$ to the groupoid of $N$-tuples of pairs $( (E_i, \psi_i) )_{i=1}^N$ such that the sum of the traces $\sum_{i=1}^N \tr(\psi_i) \in H^0(T, \cO_T)$ is zero.
\end{defn}

\begin{lemma} \label{lemma: hitchin morphism factors through Fuchs}
    The morphism $\vec{res}: \GHiggs_{g,n}^N \to\left( \mathfrak{gl}_N/\GL_N\right)^n$ factors through the closed substack $\mathscr{F}_N$.
\end{lemma}
\begin{proof}
This is essentially the residue theorem; we provide a proof in this relative setting, since we have not seen an explicit statement at this level of generality in the literature.

    It suffices to show that for any given affine scheme $T$ and any $T$-family $(C, \sigma_{\bullet}, E, \psi)$ of meromorphic Gieseker Higgs bundles, we have $\sum_{i=1}^n tr(r_i(\psi)) =0$. By taking the trace of $\psi$, we obtain a section $tr(\psi): \cO_C \to \omega_{C/T}^{log}$ whose restriction to the $i^{th}$ marking $\sigma_i$ yields $tr(r_i(\psi))$ under the residue isomorphism. Consider the short exact sequence
    \[ 0 \to \omega_{C/T} \to \omega_{C/T}^{log} \to \bigoplus_{i=1}^n (\sigma_i)_* \sigma_i^*\omega_{C/T}^{log} \to 0.\]
    Using the residue isomorphism, we may rewrite this as
    \[ 0 \to \omega_{C/T} \to \omega_{C/T}^{log} \to \bigoplus_{i=1}^n (\sigma_i)_* \cO_T \to 0.\]
    Taking global sections, we get an exact sequence
    \[H^0(C,\omega_{C/T}^{log}) \to \bigoplus_{i=1}^n H^0(T,\cO_T) \xrightarrow{\delta} H^1(C, \omega_{C/T}). \]
    By Grothendieck duality, the $\cO_T$-module $H^1(T, \omega_{C/T})$ is canonically isomorphic to $H^0(T, \cO_T)^{\vee} \cong \cO_T \cong H^0(T, \cO_T)$, and so we may rewrite the exact sequence as
    \[H^0(C,\omega_{C/T}^{log}) \to \bigoplus_{i=1}^n H^0(T,\cO_T) \xrightarrow{\delta} H^0(C,\cO_T). \]
    Under these identifications, we have that the global section $tr(\psi)$ is sent to its tuple of residues $(tr(\psi_i))_{i=1}^n$. Furthermore, we have that $\delta(tr(\psi_i))_{i=1}^n = \sum_{i=1}^n tr(\psi_i)$ by the compatibility of residues with the trace map in Grothendieck duality (see \cite[\S1]{conrad-residuevstrace}). Since the composition of the morphisms in the exact sequence is $0$, we conclude that $0 = \delta(tr(\psi_i))_{i=1}^n = \sum_{i=1}^n tr(\psi_i)$.
\end{proof}

The rest of this subsection is dedicated to proving the following.
\begin{prop} \label{prop: syntomicity of the residue morphism}
    The morphism $\vec{res}: \GHiggs_{g,n}^N \to \mathscr{F}_N \times \overline{\cM}_{g,n}$ is syntomic of pure relative dimension $N^2(2g-2+n)+1$.
\end{prop}
As a first step in the proof of \Cref{prop: syntomicity of the residue morphism}, let us prove the result over the open substack $\cM_{g,n} \subset \overline{\cM}_{g,n}$ parametrizing smooth $n$-pointed curves. For the next lemma, we denote by $\GHiggs_{g,n}^{N,\circ} \subset \GHiggs_{g,n}^N$ the open preimage of $\cM_{g,n}$ under the morphism $\GHiggs_{g,n}^N \to \overline{\cM}_{g,n}$.
\begin{lemma} \label{lemma: syntomicity of residue morphism smooth case}
    The restriction $\vec{res}: \GHiggs_{g,n}^{N,\circ} \to \mathscr{F}_N \times \cM_{g,n}$ is syntomic of pure relative dimension $N^2(2g-2+n)+1$.
\end{lemma}
\begin{proof}
    If $(C, \sigma_{\bullet})$ is a semistable $n$-pointed curve such that its stabilization $(\overline{C}, \overline{\sigma}_{\bullet})$ is smooth, then this actually forces equality $C = \overline{C}$. It follows that the stack $\GHiggs_{g,n}^{N,\circ}$ parametrizes tuples $(C, \sigma_{\bullet}, E, \psi)$, where $(C, \sigma_{\bullet})$ is a stable $n$-pointed genus $g$ smooth curve and $(E, \psi)$ is a meromorphic Higgs bundle on it. By standard deformation theory of Higgs bundles as developed in \cite{biswas-ramanana-infinitesimal-hitchin} (see also \cite[Prop. 3.8]{herrero2023meromorphic}), the forgetful morphism $\GHiggs_{g,n}^{N,\circ} \to \cM_{g,n}$ admits a relative perfect obstruction theory of virtual rank $N^2(2g-2+n)$. Since $\mathscr{F}_N \times \overline{\cM}_{g,n} \to \overline{\cM}_{g,n}$ is smooth of relative dimension $-1$, it suffices to show that $\vec{res}$ has equidimensional fibers of 
    the
    expected dimension $N^2(2g-2+n)+1$. Note that both stacks $\GHiggs_{g,n}^{N,\circ}$ and $\mathscr{F}_N$ admit an action of $\mathbb{G}_m$ that acts by scaling the Higgs field and the tuple of endomorphisms, respectively. The morphism $\vec{res}$ is $\mathbb{G}_m$-equivariant. Using the scaling action, we may include any $k$-point $p \in \GHiggs_{g,n}^{N,\circ}(k)$ corresponding to a tuple $(C, \sigma_{\bullet}, E, \psi)$ into the $\mathbb{A}^1$-family 
    \[ \mathbb{A}^1_k \to \GHiggs_{g,n}^{N,\circ}, \; \; t \mapsto (C, \sigma_{\bullet}, E, t\cdot \psi)\]
    so that the closed point $0 \in \mathbb{A}^1_k$ is mapped to the zero Higgs bundle $(C, \sigma_{\bullet}, E,0)$, which lies over the origin $\vec{0} \in \mathscr{F}_N$. By upper-semicontinuity of fiber dimension, 
we are reduced to showing that for 
    an
    $n$-pointed stable smooth curve $(C, \sigma_{\bullet})$ corresponding to a point $x \in \cM_{g,n}(k)$, the fiber of $\vec{res}$ over $\vec{0}\times x: \Spec(k) \to \mathscr{F}_N \times \cM_{g,n}$ has the expected dimension. Such fiber $\vec{res}^{-1}(\vec{0} \times x)$ parametrizes tuples $(E, \psi, u_{\bullet})$ consisting of a Higgs bundle $(E, \psi: E \to E \otimes \omega_{C})$ and a tuple of trivializations $u_i: E|_{\sigma_i} \xrightarrow{\sim} \cO_{\sigma_i}^{\oplus N}$. Therefore, $\vec{res}^{-1}(\vec{0}\times x)$ is smooth of relative dimension $N^2\cdot n$ over the stack $\Higgs(C)$ of rank $N$ Higgs bundles on the smooth curve.  The stack $\Higgs(C)$ is known to be equidimensional of dimension $N^2(2g-2)+1$, and hence it follows that $\vec{res}^{-1}(\vec{0}\times x)$ is equidimensional of 
        the
    expected dimension $N^2(2g-2+n)+1$.
\end{proof}

\begin{proof}[Proof of \Cref{prop: syntomicity of the residue morphism}]
    The morphism $\vec{res}$ factors as:
    \[\vec{res}: \GHiggs_{g,n}^N \xrightarrow{\vec{res}} \mathscr{F}_N \times \cM^{ss}_{g,n} \xrightarrow{id \times \stab} \mathscr{F}_N \times \overline{\cM}_{g,n},\]
    where $\vec{res}: \GHiggs_{g,n}^N \to \mathscr{F}_N \times \cM^{ss}_{g,n}$ sends a point $(C, \sigma_{\bullet}, E, \psi)$ to $((res_i)_{i=1}^n, (C, \sigma_{\bullet}))$. Since the stabilization morphism $\mathscr{F}_N \times \cM^{ss}_{g,n} \to \mathscr{F}_N \times \overline{\cM}_{g,n}$ is syntomic of pure relative dimension $0$, it suffices to show that $\vec{res}: \GHiggs_{g,n}^N \to \mathscr{F}_N \times \cM^{ss}_{g,n}$ is syntomic of pure relative dimension $N^2(2g-2+n)+1$. Standard deformation theory of Higgs bundles implies that the forgetful morphism $\GHiggs_{g,n}^N \to \cM^{ss}_{g,n}$ admits a relative perfect obstruction theory of virtual rank $N^2(2g-2+n)$. Since $\mathscr{F}_N$ is smooth of dimension $-1$, it suffices to show that the fibers of $\vec{res}: \GHiggs_{g,n}^N \to \mathscr{F}_N \times \cM^{ss}_{g,n}$ are equidimensional of 
        the
expected dimension $N^2(2g-2+n)+1$. Using the $\mathbb{G}_m$-action and upper-semicontinuity as in the proof of \Cref{lemma: syntomicity of residue morphism smooth case}, we may reduce to showing that for any fixed semistable $n$-pointed curve $(C, \sigma_{\bullet})$ corresponding to a point $x \in \cM^{ss}_{g,n}(k)$, the fiber $\vec{res}^{-1}(\vec{0} \times x)$ is equidimensional of expected dimension. Similarly, as in the proof of \Cref{lemma: syntomicity of residue morphism smooth case}, the fiber $\vec{res}^{-1}(\vec{0} \times x)$ parametrizes tuples $(E, \psi: E \to E \otimes \omega_{C}, u_{\bullet})$ consisting of a Higgs bundle $(E, \psi)$ on $C$ and a tuple of trivializations of the vector bundle $E$ at each point $\sigma_i$. Therefore, it suffices to show that the stack $\Higgs(C)$ parametrizing Higgs bundles $(E, \psi: E \to E \otimes \omega_{C})$ on $C$ is equidimensional of dimension $N^2(2g-2)+1$.
    
    Let $\widetilde{C} \to C$ denote the normalization of $C$, which we may write as a disjoint union of connected smooth curves $\widetilde{C} = \bigsqcup_{j \in J} C_j$. Let $n_j$ denote the number of preimages of nodes of $C$ under $C_j \to C$. We consider each $C_j$ as an $n_j$-pointed curve $(C_j, \sigma^j_{\bullet})$, where we choose an ordering of the $n_j$-preimages of nodes. Let $\Higgs(C_j, \sigma^j_{\bullet})$ denote the stack parametrizing meromorphic Higgs bundles $(E, \psi: E \to E \times \omega_{C_j}^{log})$ on the pointed smooth curve $(C_j, \sigma^j_{\bullet})$. Pulling back under the normalization morphism induces a map of stacks $\Higgs(C) \to \prod_{j \in J}\Higgs(C_j, \sigma^j_{\bullet})$. For each $j \in J$, there is a residue morphism $\vec{res}_j: \Higgs(C_j, \sigma^j_{\bullet}) \to \mathscr{F}_{N,j} \subset \left(\mathfrak{gl}_N/\GL_N\right)^{n_j}$. By the proof of \Cref{lemma: syntomicity of residue morphism smooth case}, $\vec{res}_j$ is flat of pure relative dimension $N^2(2g_j-2+n_j)+1$. We may describe the stack $\Higgs(C)$ by gluing tuples of Higgs bundles in $\prod_{j \in J} \Higgs(C_j, \sigma^j_{\bullet})$ as follows. Each node $x_i$ of $C$ corresponds to two preimages $p_i, q_i$ in $\widetilde{C} = \bigsqcup_{j \in J} C_j$. We consider the product $\left(\mathfrak{gl}_N/\GL_N\right)^2$ parametrizing pairs of residues at $p_i$ and $q_i$, and set $\Delta^{-}: \mathfrak{gl}_N/\GL_N \to \left(\mathfrak{gl}_N/\GL_N\right)^2$ to be the antidiagonal morphism given by $(E, M) \mapsto ((E,M), (E,-M))$. If we take the product as we range over all the preimages of nodes $(x_i)_{i \in I}$ of $C$, we get a morphism
    \[ \vec{\Delta}^-: \prod_{i \in I} \left(\mathfrak{gl}_N/\GL_N\right) \xrightarrow{\prod_{i}\Delta^-} \prod_{i \in I} \left(\mathfrak{gl}_N/\GL_N\right)^2 \xrightarrow{\sim}  \prod_{j \in J} \left(\mathfrak{gl}_N/\GL_N\right)^{n_j},\]
    where in the target the product ranges over all marked points in every $C_j$. We denote by $\mathscr{Z} \subset \prod_{i \in I} \left(\mathfrak{gl}_N/\GL_N\right)$ the preimage under $\vec{\Delta}^-$ of the closed substack $\prod_{j \in J} \mathscr{F}_{N,j} \subset \prod_{j \in J} \left(\mathfrak{gl}_N/\GL_N\right)^{n_j}$, which imposes that the sum of the traces of the residues over each $C_j$ is zero. Then, because every Higgs bundle in $\Higgs(C)$ can be obtained by gluing Higgs bundles over the normalization $\widetilde{C}$, it follows that $\Higgs(C)$ fits into a Cartesian diagram:
    \[
	\begin{tikzcd}
		\Higgs(C) \ar[r] \ar[d] & \prod_{j \in J} \Higgs(C_j, \sigma^j_{\bullet}) \ar[d, "\prod_{j} \vec{res}_j"] \\ \mathscr{Z}  \ar[r, "\vec{\Delta}^-"] & \prod_{j \in J}\mathscr{F}_{N,j}
	\end{tikzcd}
\]
    Since each $\vec{res}_j$ is flat of pure relative dimension $N^2(2g_j-2) + 1$, it follows that $\Higgs(C) \to \mathscr{Z}$ is flat or pure relative dimension $N^2 (2\sum_{j \in J}g_j -2|J| + 2|I|) + |J|$. On the other hand, $\mathscr{Z} \subset \prod_{i \in I} \left(\mathfrak{gl}_N/\GL_N\right)$ is cut out by some linear equations (sum of certain traces) in $\prod_{i \in I} \mathfrak{gl}_N$, and hence it is smooth (equidimensional) of some given dimension $\dim(\mathscr{Z})$. Each of the $|J|$-components $C_j$ induces a linear equation corresponding to $\mathscr{F}_{N,j} \subset \left(\mathfrak{gl}_N/\GL_N\right)^{n_j}$, but there is one linear relation among such equations, and therefore the codimension of $\mathscr{Z} \subset \prod_{i \in I} \left(\mathfrak{gl}_N/\GL_N\right)^{n_j}$ is $|J|-1$. We conclude that $\dim(\mathscr{Z}) = 1-|J|$, and hence $\dim(\Higgs(C)) = N^2(2(\sum_{j \in J} g_j -|J| + |I| + 1) - 2) +1$. Since the genus $g$ of $C$ is given by $g=\sum_{j \in J} g_j -|J| + |I| + 1$, we conclude that $\Higgs(C)$ is equidimensional of dimension $N^2(2g-2)+1$, as desired.
\end{proof}

\begin{subsection}{Modified Hitchin base for fixed nilpotent residues}
Let $\cO_{\bullet}= (\cO_i)_{i=1}^n$ be a tuple of nilpotent conjugacy classes. In this subsection, we will explain a version $A^{\cO_{\bullet}}$ of the Hitchin base developed in \cite{BDD-modified-hitchin-gln}, which will serve as the target of a ``log-symplectic system'' for the stack $\GHiggs_{g,n}^{\cO_{\bullet}}$ of meromorphic Gieseker Higgs bundles with fixed residue classes $\cO_{\bullet}$ (as in the upcoming \Cref{defn: gieseker higgs bundles with fixed residues}). Recall that any nilpotent conjugacy class in $\mathfrak{gl}_N$ corresponds naturally via its Jordan block decomposition to a partition $\vec{\lambda}$ of $N$, where $\vec{\lambda}= (N_1, N_2, \ldots, N_l)$ is a tuple of positive integers with $N_1 \geq N_2 \geq \ldots \geq N_l$ and $\sum_{i=1}^l N_i =N$.

\begin{defn}[Vanishing and pole orders]
    Let $\vec{\lambda}= (N_1, N_2, \ldots, N_l)$ be a partition of $N$, with $N_1 \geq N_2 \geq $. For any $1 \leq j \leq N$, we define $\chi^{\vec{\lambda}}_j$ to be the unique index such that 
    \[\sum_{i=1}^{\chi^{\vec{\lambda}}_j-1} N_i < j \leq \sum_{i=1}^{\chi^{\vec{\lambda}}_j} N_i \]
    We set $\pi^{\vec{\lambda}}_j:= N- \chi^{\vec{\lambda}}_i$.
\end{defn}

Fix a tuple $\cO_{\bullet}= (\cO_i)_{i=1}^n$ of nilpotent conjugacy classes. We denote by $\vec{\lambda}_i$ the partition of $N$ corresponding to the Jordan block decomposition of $\cO_i$.

\begin{defn}[Naive Hitchin base]
    The naive Hitchin base $\cA^{\cO_{\bullet}}_{nv}$ with fixed residues $\cO_{\bullet}$ is the subfunctor of $\cA$ that sends a scheme $T$ into the groupoid of tuples $(C, \sigma_{\bullet}, (a_j)_{j=1}^N)$ such that 
\[a_j \in H^0\left(C, \omega_{C/T}^{\otimes j}(\sum_{i=1}^n \pi^{\vec{\lambda}_i}_{j} \sigma_i)\right) \subset H^0\left(C, (\omega_{C/T}^{log})^{\otimes j}\right).\]
\end{defn}
The naive Hitchin base $\cA^{\cO_{\bullet}}_{nv}$ is represented by a stack that is affine and of finite type over $\overline{\cM}_{g,n}$. However, unlike in the case of $\cA$, the stack $\cA^{\cO_{\bullet}}_{nv}$ is not the total space of a vector bundle over $\overline{\cM}_{g,n}$. This phenomenon appears in the case when $g=0$ because the fibers of the vector spaces $H^0\left(C, \omega_{C/T}^{\otimes j}(\sum_{i=1}^n \pi^{\vec{\lambda}_i}_{j} \sigma_i)\right)$ might jump as we vary the $n$-pointed curve $(C, \sigma_{\bullet})$; see \cite[\S3.3]{BDD-modified-hitchin-gln} for an example in the case of $\overline{\cM}_{0,4}$.

In order to remedy this jumping phenomenon, a modified Hitchin base $\cA^{\cO_{\bullet}}$ is proposed in \cite{BDD-modified-hitchin-gln}. In order to recall it, we need to set up some notation first. Let $\cM_{g,n} \subset \overline{\cM}_{g,n}$ the open substack that parametrizes smooth stable $n$-pointed curves. The closed complement $\Delta = \overline{\cM}_{g,n} \setminus \cM_{g,n}$ is a normal crossings divisor on $\overline{\cM}_{g,n}$. We have a decomposition into irreducible components $\Delta = \Delta_{irr} \cup \bigcup_{h \in H} \Delta_h$, where the general point of $\Delta_{irr}$ parametrizes irreducible curves with a single node, and where the index $h$ runs over the set $H$ of tuples $(g_1, g_2, S_1, S_2)$ 
consisting of the following:
\begin{enumerate}[(1)]
    \item A pair of integers $g_1, g_2 \geq 0$ such that $g_1 +g_2 =g$.
    \item A pair of subsets $S_1, S_2 \subset \{1, 2, \ldots , n\}$ such that $S_1 \sqcup S_2 = \{1,2,\ldots , n\}$ and such that $|S_i| \geq 2$ if $g_i =0$.
\end{enumerate}
The general point of $\Delta_h$ parametrizes stable $n$-pointed curves $(C, \sigma_{\bullet})$ where the curve $C$ is the union of two irreducible smooth subcurves $C = C_1 \cup C_2$ meeting at one node, and we have equalities $g(C_i)=g_i$ 
for the genera and the property that $S_j = \{ i \, \mid \, \sigma_i \in C_j \}$. 

\begin{remark} \label{remark: redundancies indexes}
    We remark that there are some redundancies in our indexing of the irreducible components $\Delta_h$ above. Indeed, if $g_1=g_2= \frac{1}{2}g$, then the tuples $(\frac{1}{2}g, \frac{1}{2}g, S_1, S_2)$ and $(\frac{1}{2}g, \frac{1}{2}g, S_2, S_1)$ yield the same irreducible component.
\end{remark}

We will be interested in components $\Delta_h$ obtained as follows. For any given subset $S \subset \{1, 2, \ldots, n\}$ of cardinality at least $2$, we set $S^{\vee}:= \{1,2,\ldots, n\} \setminus S$ and set $h_S= (0,g,S, S^{\vee})$. For any such $S$, the fiber $\cC_{\eta}$ over the generic point $\eta \in \Delta_{h_S}$ has a unique genus $0$ component containing exactly the markings indexed by $i \in S$. We denote by $\cC_S \subset \cC$ the divisor on the (smooth) total space of the universal curve $\cC$ obtained by taking the closure of such component in $\cC_{\eta} \hookrightarrow \cC$. We denote by $\cO(\cC_S)$ the corresponding line bundle on $\cC$.

\begin{notn}[Some relevant divisors on $\cC$]
    Fix $S \subset \{1,2,\ldots, n\}$ of cardinality at least $2$. Let $T \to \overline{\cM}_{g,n}$ be a morphism classifying an $n$-pointed stable curve $(C, \overline{\sigma}_{\bullet})$, which we may think as the pullback $(\cC_T, (\sigma_{\bullet})_T)$ over $T$ of the universal family $(\cC, \sigma_{\bullet}) \to \overline{\cM}_{g,n}$. We set $\cO(C_S)$ to denote the pullback of the line bundle $\cO(\cC_S)$ under the morphism $C = \cC_T \to \cC$.
\end{notn}

\begin{defn}[Twist orders]
    Fix an $n$-tuple of nilpotent conjugacy classes $\cO_{\bullet}$ with corresponding tuple of partitions $(\vec{\lambda}_{\bullet})$. For any given $S \subset \{1,2,\ldots, n\}$ of cardinality at least $2$, and any given $1 \leq j \leq N$, we set
    \[ n^S_{\vec{\lambda}_{\bullet},j}= max\left(0, j-1-\sum_{i \in S} (j-\chi_j^{\vec{\lambda}_i})\right)\]
\end{defn}

Now we are ready to define the modified Hitchin base as in \cite{BDD-modified-hitchin-gln}.
\begin{defn}[Modified Hitchin base]
    Let $\cO_{\bullet}= (\cO_i)_{i=1}^n$ be a tuple of nilpotent conjugacy classes, with corresponding Jordan partition $\vec{\lambda}_i$. The Hitchin base $\cA^{\cO_{\bullet}}$ with fixed residue classes $\cO_{\bullet}$ is the pseudofunctor from schemes into groupoids that sends a scheme $T$ to the groupoid of tuples $(C, \sigma_{\bullet}, (a_j)_{j=1}^N)$, where $(C, \sigma_{\bullet})$ is a $T$-family of $n$-pointed stable curves in $\overline{\cM}_{g,n}(T)$ and $a_j$ is a global section
    \[ a_j: \cO_{C} \to \omega_{C/T}^{\otimes j}(\sum_{i=1}^n \pi^{\vec{\lambda}_i}_{j} \sigma_i) \otimes \bigotimes_{S} \cO(C_S)^{\otimes{-n^S_{ \vec{\lambda}_{\bullet}, j}}},\]
  where the index $S$ runs over all subsets $S \subset \{1, 2, \ldots, n\}$ of cardinality at least $2$.
\end{defn}

The following is one of the main results in \cite{BDD-modified-hitchin-gln}.
\begin{thm}
    Let $\cO_{\bullet}= (\cO_i)_{i=1}^n$ be a tuple of nilpotent conjugacy classes, with corresponding Jordan partition $\vec{\lambda}_i$. The forgetful morphism
    \[ \cA^{\cO_{\bullet}} \to \overline{\cM}_{g,n}, \; \; \; (C, \sigma_{\bullet}, a_{\bullet}) \mapsto (C, \sigma_{\bullet})\]
    exhibits $\cA^{\cO_{\bullet}}$ as the total space of a vector bundle of rank $N^2(g-1) + 1+  \frac{1}{2}\sum_{i=1}^n\dim(\cO_i)$ on $\overline{\cM}_{g,n}$.
\end{thm}

By construction, there is a natural morphism $\cA^{\cO_{\bullet}} \to \cA^{\cO_{\bullet}}_{nv}$ induced by the following inclusions of line bundles on the universal curve $(\cC, \sigma_{\bullet}) \to \overline{\cM}_{g,n}$
\[ \omega_{\cC/\overline{\cM}_{g,n}}^{\otimes j}(\sum_{i=1}^n \pi^{\vec{\lambda}_i}_{j} \sigma_i) \otimes \bigotimes_S \cO(\cC_S)^{\otimes{-n^S_{ \vec{\lambda}, j}}} \hookrightarrow \omega_{\cC/\overline{\cM}_{g,n}}^{\otimes j}(\sum_{i=1}^n \pi^{\vec{\lambda}_i}_{j} \sigma_i)\]
for any given index $1 \leq j \leq N$.
The following observation will be useful for us in defining the Hitchin morphism with fixed residues.

\begin{lemma} \label{lemma: A naive vs A modified}
    Let $T \to \overline{\cM}_{g,n}$ be a flat scheme over $\overline{\cM}_{g,n}$. Then the morphism of affine $\overline{\cM}_{g,n}$-stacks $\cA^{\cO_{\bullet}} \to \cA^{\cO_{\bullet}}_{nv}$ induces a bijection on the sets of sections $\cA^{\cO_{\bullet}}(T) \xrightarrow{\sim} \cA^{\cO_{\bullet}}_{nv}(T)$.
\end{lemma}
\begin{proof}
    To simplify notation, for each $1 \leq j \leq N$ we use the notation 
    \[\cL_j:= \omega_{\cC/\overline{\cM}_{g,n}}^{\otimes j}(\sum_{i=1}^n \pi^{\vec{\lambda}_i}_{j} \sigma_i), \; \; \text{ and }\;\; \cL'_j := \omega_{\cC/\overline{\cM}_{g,n}}^{\otimes j}(\sum_{i=1}^n \pi^{\vec{\lambda}_i}_{j} \sigma_i) \otimes \bigotimes_S \cO(\cC_S)^{\otimes{-n^S_{ \vec{\lambda}, j}}}.\]
    The functor $\cA^{\cO_{\bullet}}$ classifies sections of (pullbacks of) $\bigoplus_{j=1}^N \cL'_j$, whereas $\cA^{\cO_{\bullet}}_{nv}$ classifies  sections of (pullbacks of) $\bigoplus_{j=1}^N \cL_j$. The morphism $\cA^{\cO_{\bullet}} \hookrightarrow \cA^{\cO_{\bullet}}_{nv}$ is induced by the injection of sheaves $\bigoplus_{j=1}^N \cL'_j \hookrightarrow \bigoplus_{j=1}^N\cL_j$ on $\cC$. Let $\pi: \cC \to \overline{\cM}_{g,n}$ denote the structure morphism of $\cC$. By the compatibility of the formation of the pushforwards $\pi_*(\bigoplus_{j=1}^N \cL'_j)$ and $\pi_*(\bigoplus_{j=1}^N \cL_j)$ with flat base-change, in order to prove the lemma for any flat $T \to \overline{\cM}_{g,n}$ it suffices to show that the corresponding injections of sheaves $\varphi_j: \pi_*(\cL'_j) \hookrightarrow \pi_*(\cL_j)$ are isomorphisms over $\overline{\cM}_{g,n}$.

    Notice that $\pi_*(\cL'_j)$ is a vector bundle on $\overline{\cM}_{g,n}$ by the main theorem of \cite{BDD-modified-hitchin-gln}, whereas $\pi_*(\cL_j)$ is a torsion-free sheaf on $\overline{\cM}_{g,n}$ by definition. By Hartogs'
 theorem it suffices to check surjectivity of $\varphi_j$ away from a codimension 2 substack. The morphism $\varphi_j$ is an isomorphism over $\overline{\cM}_{g,n} \setminus \bigcup_{S} \Delta_{h_S}$, 
 where $S$ runs over all subsets of $\{1, 2, \ldots, n\}$ of cardinality at least $2$ (see the paragraph after \Cref{remark: redundancies indexes} for the definition of $h_S$). Therefore, it suffices to show that for any given $S$, the morphism $\varphi_j$ is an isomorphism at sufficiently general points over the divisor $\Delta_{h_S}$. By smoothness of $\overline{\cM}_{g,n}$, over sufficiently general points of $\Delta_{h_S}$, both torsion-free sheaves $\pi_*(\cL'_j)$ and $\pi_*(\cL_j)$ are vector bundles of the same rank. 
 For any given $S$ and any sufficiently general point $x: \Spec(k) \to \Delta_{h_S}$, to show the required isomorphism of vector bundles locally around the image of $x$ we are reduced to showing that the induced morphism on fibers $(\varphi_j)_x: \pi_*(\cL'_j)_x \to \pi_*(\cL_j)_x$ is injective. 
    
    Let us we denote by $\cC_x$ the fiber of $\cC \to \overline{\cM}_{g,n}$. By the main theorem in \cite{BDD-modified-hitchin-gln}, the formation of the pushforward $\pi_*(\cL'_j)$ commutes with arbitrary base-change, and hence $\pi_*(\cL'_j)_x = H^0(\cC_x,\cL'_j|_{\cC_x})$. Using this, we see that it suffices to show that the induced morphism $\psi: H^0(\cC_x,\cL'_j|_{\cC_x}) \to H^0(\cC_x,\cL_j|_{\cC_x})$ is injective. If $x$ does not lie on $\Delta_{h_{S'}}$ for $S' \neq S$, then we have $\cL'_j|_{\cC_x} = \cL_j(-n^S_{\vec{\lambda}_{\bullet},j}\cC_S)|_{\cC_x}$. We may factor the morphism $\psi$ into a chain of morphisms $\psi = \psi_1 \circ \psi_2 \circ \ldots \circ \psi_{n^S_{\vec{\lambda}_{\bullet},j}}$ 
    as follows:
    \[ \psi_{n^S_{\vec{\lambda}_{\bullet},j}}: H^0(\cC_x,\cL_j'|_{\cC_x}) = H^0(\cC_x,\cL_j(-n^S_{\vec{\lambda}_{\bullet},j}\cC_S)|_{\cC_x}) \to  H^0(\cC_x,\cL_j(-(n^S_{\vec{\lambda}_{\bullet},j}-1)\cC_S)|_{\cC_x})\]
    \[\psi_{n^S_{\vec{\lambda}_{\bullet},j}-1}: H^0(\cC_x,\cL_j(-(n^S_{\vec{\lambda}_{\bullet},j}-1)\cC_S)|_{\cC_x}) \to  H^0(\cC_x,\cL_j(-(n^S_{\vec{\lambda}_{\bullet},j}-2)\cC_S)|_{\cC_x})\]
    \[ \ldots\]
    \[ \psi_1: H^0(\cC_x,\cL_j(-\cC_S)|_{\cC_x}) \xrightarrow{\psi_1} H^0(\cC_x,\cL_j|_{\cC_x}).\]
    We shall conclude the proof by showing that $\psi_i$ is injective for all $i$. Recall that the sufficiently general point $x: \Spec(k) \to \Delta_S$ parametrizes a curve $\cC_x$ with two irreducible components $\cC_x = C_1 \cup C_2$ meeting at one node. Here $C_1 = \cC_x \cap \cC_S$
    is a smooth curve of genus $0$. The line bundle restriction $\cO(\cC_S)|_{C_1}$ has degree $-1$, whereas the restriction $\cO(\cC_S)|_{C_2}$ has degree $1$. For any given number $1 \leq i \leq n^S_{\vec{\lambda}_{\bullet},j}$, the kernel of the morphism
    \[\psi_i: H^0(\cC_x,\cL_j(-i\cC_S)|_{\cC_x}) \to H^0(\cC_x, \cL_j(-(i-1)\cC_S)|_{\cC_x})\]
    consists of all the global sections of the line bundle $\cL_j(-i\cC_S)|_{\cC_x}$ that are completely supported on $C_1$. Therefore, it suffices to show the following:
    
    \noindent \textbf{Claim:} For any index $1 \leq i \leq n^S_{\vec{\lambda}_{\bullet},j}$, there are no nonzero global sections inside $H^0(\cC_x, \cL_j(-(i-1)\cC_S)|_{\cC_x})$ which are completely supported on the component $C_1$. 
    
    Since the line bundle $\cO(\cC_S)$ is negative on $C_1$, if the claim follows for an index $i$, then it also follows for the smaller indexes. Therefore, it suffices to show the claim when $i=n^S_{\vec{\lambda}_{\bullet},j}$. This amounts to showing that there are no nonzero global sections of $\cL'_j|_{\cC_x}$ that are supported on $C_1$, which follows directly from \cite[Equation (127), pg.67]{BDD-modified-hitchin-gln} (In fact Equation (127) in the reference states the a priori stronger statement that all sections of $\cL'_j|_{\cC_x}$ are completely supported on the complementary irreducible component $C_2$).
\end{proof}

\end{subsection}

\begin{subsection}{Hitchin morphism for fixed nilpotent residues} \label{subsection: hitchin morphism for fixed nilpotent residues}
    Let $\cO_{\bullet} = (\cO_i)_{i=1}^n$ be an $n$-tuple of conjugacy classes in $\mathfrak{gl}_N$. We have a locally closed immersion 
    \[\vec{\iota}: \prod_{i=1}^N (\cO_i/\GL_N) \hookrightarrow (\mathfrak{gl}_N/\GL_N)^n.\]
    We suppose that we have $\sum_{i=1}^N tr(\cO_i)=0$, so that $\vec{\iota}$ factors through the Fuchs substack $\mathscr{F}_N \subset \prod_{i=1}^n (\mathfrak{gl}_N/\GL_N)$.

    \begin{defn}[Stack of meromorphic Gieseker Higgs bundles with fixed residues] \label{defn: gieseker higgs bundles with fixed residues}
        Let $\cO_{\bullet}$ be a tuple of conjugacy classes with $\sum_{i=1}^n \tr(\cO_i) =0$. We denote by $\GHiggs_{g,n}^{\cO_{\bullet}}$ the stack fitting into the following Cartesian diagram
          \[
	\begin{tikzcd}
		\GHiggs_{g,n}^{\cO_{\bullet}} \ar[r] \ar[d] & \GHiggs_{g,n}^N \ar[d, "\vec{res}"] \\ (\prod_{i=1}^N (\cO_i/\GL_N)) \times \overline{\cM}_{g,n} \ar[r, "\vec{\iota} \times id"] & \mathscr{F}_N \times \overline{\cM}_{g,n}
	\end{tikzcd}
\]
It is the locally closed substack of $\GHiggs_{g,n}^N$ that classifies Gieseker Higgs bundles such that the $i^{th}$ residue of the Higgs field lies in the conjugacy class $\cO_i$.
    \end{defn}

    We shall also employ the closure of such locally closed substacks. These can be obtained as fiber products by using the closures $\overline{\cO}_i \subset \mathfrak{gl}_N$ of the conjugacy classes.
    \begin{defn}[Closure stack]
        Let $\cO_{\bullet}$ be a tuple of conjugacy classes satisfying $\sum_{i=1}^n \tr(\cO_i)=0$. We denote by $\overline{\GHiggs}_{g,n}^{\cO_{\bullet}}$  the closed substack of $\GHiggs_{g,n}^N$ fitting into the following Cartesian diagram
          \[
	\begin{tikzcd}
		\overline{\GHiggs}_{g,n}^{\cO_{\bullet}} \ar[r] \ar[d] & \GHiggs_{g,n}^N \ar[d, "\vec{res}"] \\ (\prod_{i=1}^N (\overline{\cO}_i/\GL_N)) \times \overline{\cM}_{g,n} \ar[r, "\vec{\iota} \times id"] & \mathscr{F}_N \times \overline{\cM}_{g,n}
	\end{tikzcd}
\]
    \end{defn}

\begin{lemma} \label{lemma: syntomicity of gieseker higgs bundles with fixed residues}
    The structure morphism $\overline{\GHiggs}_{g,n}^{\cO_{\bullet}} \to \overline{\cM}_{g,n}$ is flat with Cohen-Macaulay fibers of relative dimension $N^2(2g-2)+1 + \sum_{i=1}^N\dim(\cO_i)$. 
    
    Furthermore, the restriction of this morphism to the open substack $\GHiggs_{g,n}^{\cO_{\bullet}} \subset \overline{\GHiggs}_{g,n}^{\cO_{\bullet}}$ is syntomic with equidimensional fibers of the same dimension $N^2(2g-2)+1 + \sum_{i=1}^N\dim(\cO_i)$.
\end{lemma}
\begin{proof}
    Given the definition of $\overline{\GHiggs}_{g,n}^{\cO_{\bullet}}$
    as a fiber product, the statement of this lemma for $\overline{\GHiggs}_{g,n}^{\cO_{\bullet}} \to \overline{\cM}_{g,n}$ is a direct consequence of the fact that $\vec{res}: \GHiggs_{g,n}^N \to \mathscr{F}_n \times \overline{\cM}_{g,n}$ is syntomic of relative dimension $N^2(2g-2+n)+1$ (\Cref{prop: syntomicity of the residue morphism}) and the fact that the closed substack $\prod_{i=1}^n (\overline{\cO}_i/\GL_N)$ inside $\mathscr{F}_N$ are Cohen-Macaulay of dimension $\sum_{i=1}^n(\dim(\cO_i) - N^2)$ (indeed, each orbit closure $\overline{\cO}_i$ is Cohen-Macaulay by the main theorem in \cite{kraft-procesi}). The statement for the open substack $\GHiggs_{g,n}^{\cO_{\bullet}}$ follows similarly, using that $\prod_{i=1}^n (\cO_i/\GL_N)$ is smooth (and so in particular syntomic).
\end{proof}

\begin{prop} \label{prop: factorization of the hitchin morphism with fixed nilpotent residues}
    Suppose that $\cO_{\bullet}$ as above is a tuple of nilpotent conjugacy classes. The composition $\overline{\GHiggs}_{g,n}^{\cO_{\bullet}} \hookrightarrow \GHiggs_{g,n} \xrightarrow{H} \cA$ factors uniquely through $\cA^{\cO_{\bullet}} \to \cA$. In particular, the same holds for the composition $\GHiggs_{g,n}^{\cO_{\bullet}} \xrightarrow{H} \cA$.
\end{prop}
\begin{proof}
    By \Cref{lemma: syntomicity of gieseker higgs bundles with fixed residues}, the morphism $\overline{\GHiggs}_{g,n}^{\cO_{\bullet}} \to \overline{\cM}_{g,n}$ is flat. Therefore, by \Cref{lemma: A naive vs A modified}, it suffices to show that the morphism $\overline{\GHiggs}_{g,n}^{\cO_{\bullet}} \xrightarrow{H} \cA$ factors through the closed subfunctor $\cA^{\cO_{\bullet}}_{nv}$. This amounts to a local computation around the markings to ensure that having a tuple of residues in the closures of the conjugacy classes of the tuple $\cO_{\bullet}$ implies the imposed vanishing orders in the definition of $\cA^{\cO_{\bullet}}_{nv}$ as explained in \cite[\S2]{BDD-modified-hitchin-gln}, see \cite[Thm. 7]{baraglia-kamgarpour} for an explicit computation. 
\end{proof}

\end{subsection}

\begin{subsection}{Flatness of Hitchin morphisms} \label{subsection: flatness of the hitchin morphism}
\begin{prop} \label{prop: syntomicity of the hitchin morphism}
    The morphism $\overline{\GHiggs}_{g,n}^{\cO_{\bullet}} \to \cA^{\cO_{\bullet}}$ is flat of pure relative dimension $N^2(g-1)+\frac{1}{2}\sum_{i=1}^n\dim(\cO_i)$. 
    
    Furthermore, the morphism $\GHiggs_{g,n}^{\cO_{\bullet}} \to \cA^{\cO_{\bullet}}$ is syntomic of pure relative dimension $N^2(g-1)+\frac{1}{2}\sum_{i=1}^n\dim(\cO_i)$.
\end{prop}
\begin{proof}
    The source $\overline{\GHiggs}_{g,n}^{\cO_{\bullet}}$ is flat and Cohen-Macaulay of pure relative dimension $N^2(2g-2) + 1 \sum_{i=1}^N\dim(\cO_i)$ over $\overline{\cM}_{g,n}$, whereas the base $\cA^{\cO_{\bullet}}$ is smooth and of pure relative dimension $N^2(g-1) + 1+  \frac{1}{2}\sum_{i=1}^n\dim(\cO_i)$ over $\overline{\cM}_{g,n}$. By miracle flatness \cite[\href{https://stacks.math.columbia.edu/tag/00R4}{Tag 00R4}]{stacks-project}
    , suffices to show that the fibers are equidimensional of the expected dimension $N^2(g-1)+\frac{1}{2}\sum_{i=1}^n\dim(\cO_i)$. This will also imply the syntomicity of the restriction $\GHiggs_{g,n}^{\cO_{\bullet}} \to \cA^{\cO_{\bullet}}$ to the open substack $\GHiggs_{g,n}^{\cO_{\bullet}} \subset \overline{\GHiggs}_{g,n}^{\cO_{\bullet}}$, because $\GHiggs_{g,n}^{\cO_{\bullet}}$ is syntomic over $\overline{\cM}_{g,n}$.
    
    Notice that $\mathbb{G}_m$ acts on $\overline{\GHiggs}_{g,n}^{\cO_{\bullet}}$ by scaling the Higgs field. Any $k$-point $(C, \sigma_{\bullet}, E, \psi)$ of $\overline{\GHiggs}_{g,n}^{\cO_{\bullet}}$ fits into a $\mathbb{G}_m$-equivariant $\mathbb{A}^1_k$-family 
    \[ \mathbb{A}^1_k \to \overline{\GHiggs}_{g,n}^{\cO_{\bullet}}, \; \; t \mapsto (C, \sigma_{\bullet}, E, t\psi)\]
    such that the image of $0$ lies over the zero section $0_{\overline{\cM}_{g,n}} \subset \cA^{\cO_{\bullet}}$ of the total space of the vector bundle $\cA^{\cO_{\bullet}}$. By upper-semicontinuity of fiber dimension, we are reduced to showing that the preimage $H^{-1}(0_{\overline{\cM}_{g,n}})$ has fibers over $\overline{\cM}_{g,n}$ that are equidimensional of the expected dimension $N^2(g-1)+\frac{1}{2}\sum_{i=1}^n\dim(\cO_i)$. We have a factorization
    \[H^{-1}(0_{\overline{\cM}_{g,n}})\to \cM^{ss}_{g,n} \xrightarrow{\stab} \overline{\cM}_{g,n} \]
    where $\stab$ is syntomic of pure relative dimension $0$. It suffices to show that the fibers of $H^{-1}(0_{\overline{\cM}_{g,n}}) \to \cM^{ss}_{g,n}$ are equidimensional of the expected dimension. For any given semistable curve $(C, \sigma_{\bullet})$ corresponding to a point $x \in \cM^{ss}_{g,n}$, the fiber $\left(H^{-1}(0_{\overline{\cM}_{g,n}})\right)_x$ is an open substack of the stack $\Higgs(C, \sigma_{\bullet})^{nilp, \cO_{\bullet}}$ that parametrizes meromorphic Higgs bundles $(E, \psi: E \to E \otimes \omega_{C}^{log})$ on $C$ such that $\psi$ is everywhere nilpotent and such that the $i^{th}$ residue of $\psi$ lies in $\overline{\cO}_{i}$. We are left to show that $\Higgs(C, \sigma_{\bullet})^{nilp, \cO_{\bullet}}$ is equidimensional of dimension $N^2(g-1)+\frac{1}{2}\sum_{i=1}^n\dim(\cO_i)$. (In fact, from our reductions above using upper-semicontinuity, we know that each irreducible component of $\Higgs(C, \sigma_{\bullet})^{nilp, \cO_{\bullet}}$ has dimension at least $N^2(g-1)+\frac{1}{2}\sum_{i=1}^n\dim(\cO_i)$, and hence it suffices to show that each irreducible component of $\Higgs(C, \sigma_{\bullet})^{nilp, \cO_{\bullet}}$ has dimension at most $N^2(g-1)+\frac{1}{2}\sum_{i=1}^n\dim(\cO_i)$).

    To show this, we may use a similar gluing argument as in the proof of \Cref{prop: syntomicity of the residue morphism}. We first prove the smooth case.
    
    \noindent \textbf{Claim:} Let $(C, \sigma_{\bullet})$ be an $n$-pointed smooth curve of genus $g$. Then, each irreducible component of the stack $\Higgs(C, \sigma_{\bullet})^{nilp, \cO_{\bullet}}$ has dimension at most $N^2(g-1)+\frac{1}{2}\sum_{i=1}^n\dim(\cO_i)$.

    \noindent \textit{Proof of the claim.} For each index $i$, let $\vec{\lambda}_i$ denote the partition of $N$ corresponding to the Jordan block decomposition of any matrix in the orbit $\cO_i$. Set $P_{\vec{\lambda}_i} \subset \GL_N$ to be the standard parabolic subgroup of $\GL_N$ with Levi subgroup of block diagonal matrices of sizes determine by $\vec{\lambda}_i$. We denote by $\Higgs(C,\sigma_{\bullet})_{P_{\vec{\lambda}_i}}$ the stack of (strongly) quasi-parabolic Higgs bundles with parabolic reduction to $P_{\vec{\lambda}_i}$ at the point $\sigma_i$ (as defined in \cite[Defn. 4]{baraglia-kamgarpour-varma} for the case of a single point $\sigma$). We denote by $\Higgs(C,\sigma_{\bullet})^{nilp}_{P_{\vec{\lambda}_i}} \subset \Higgs(C,\sigma_{\bullet})_{P_{\vec{\lambda}_i}}$ the closed substack of nilpotent quasi-parabolic Higgs bundles. There is a surjective morphism $\Higgs(C,\sigma_{\bullet})^{nilp}_{P_{\vec{\lambda}_i}} \to \Higgs(C, \sigma_{\bullet})^{nilp, \cO_{\bullet}}$ given by forgetting the parabolic reductions. Therefore, it suffices to show that the dimension of each irreducible component of $\Higgs(C,\sigma_{\bullet})^{nilp}_{P_{\vec{\lambda}_i}}$ is at most the expected dimension $N^2(g-1)+\frac{1}{2}\sum_{i=1}^n\dim(\cO_i)$. This is a direct consequence of \cite[Thm. 10(i)]{baraglia-kamgarpour-varma}, thus concluding the proof of the \textbf{Claim}.
    
    By our discussion above, the claim implies in particular that $\Higgs(C, \sigma_{\bullet})^{nilp, \cO_{\bullet}}$ is equidimensional of dimension $N^2(g-1)+\frac{1}{2}\sum_{i=1}^n\dim(\cO_i)$ whenever $C$ is a smooth curve. To conclude the proof of the proposition, let us address the more general case of a semistable $n$-pointed curve $(C, \sigma_{\bullet})$. 
    
    We use a similar setup as in the proof of \Cref{prop: syntomicity of the residue morphism}. We let $\widetilde{C} \to C$ be the normalization with set of components $\widetilde{C} = \sqcup_{j \in J} C_j$. We denote by $\{\widetilde{\nu}_{\bullet}^j\}$ the set of preimages of the nodes of $C$ in $C_j$, and we denote by $n_j$ the number of such preimages. For each $j \in J$, we let $(\widetilde{\sigma}^j_{l})_{ l \in L_j}$ be the preimages in $C_j$ of the markings $\sigma_{\bullet}$, where $L_j \subset \{1, \ldots, n\}$ and $\widetilde{\sigma}^j_{l}$ maps to $\sigma_l$.  We view $C_j$ as a $(|L_j|+n_j)$-pointed curve with set of markings $\{\widetilde{\nu}_{\bullet}^j \cup \widetilde{\sigma}^j_{\bullet}\}$. We denote by $\Higgs(C_j, \widetilde{\nu}^j_{\bullet}, \widetilde{\sigma}^j_{\bullet})^{nilp, \cO_{\bullet}}$ the stack parametrizing nilpotent meromorphic Higgs bundles $(E, \psi: E \to E \otimes \omega_C)$ such that for all $i \in L_j$ the residue at $\widetilde{\sigma}^j_i$ belongs to the closure $\overline{\cO}_i$. Let $\mathfrak{gl}_N^{nilp} \subset \mathfrak{gl}_N$ denote the closed subscheme of nilpotent matrices,
    and set $\mathscr{F}^{nilp}_{N,j}: = \left(\mathfrak{gl}_N^{nilp}/\GL_N\right)^{n_j}$. Then for each $j \in J$ we have a residue morphism for the $\widetilde{\nu}^j_{\bullet}$-markings $\vec{res}_j: \Higgs(C_j, \widetilde{\nu}^j_{\bullet}, \widetilde{\sigma}^j_{\bullet})^{nilp, \cO_{\bullet}} \to \mathscr{F}^{nilp}_{N,j}$. Taking the product, we obtain a morphism 
\[\prod_{j\in J} \vec{res}_j: \prod_{j \in J} \Higgs(C_j, \widetilde{\nu}^j_{\bullet}, \widetilde{\sigma}^j_{\bullet})^{nilp, \cO_{\bullet}} \to \prod_{j \in J} \left(\mathfrak{gl}^{nilp}_N/\GL_N\right)^{n_j}\]
Let $(x_i)_{i \in I}$ denote the nodes of $C$, which we index by a finite set $I$. Just as in the proof of \Cref{prop: syntomicity of the residue morphism}, as we range over the nodes $x_i$, we get a product of antidiagonal morphisms 
\[\vec{\Delta}^-: \prod_{i \in I}(\mathfrak{gl}^{nilp}/\GL_N) \to \prod_{i \in I} \left(\mathfrak{gl}^{nilp}/\GL_N\right)^2 \cong \prod_{j \in J} \left(\mathfrak{gl}^{nilp}/\GL_N\right)^{n_j}\]
    Since meromorphic Higgs bundles on $C$ arise from gluing meromorphic Higgs bundles on $\widetilde{C}$, we may fit $\Higgs(C, \sigma_{\bullet})^{nilp, \cO_{\bullet}}$ into the following Cartesian diagram:
        \[
	\begin{tikzcd}
		\Higgs(C, \sigma_{\bullet})^{nilp, \cO_{\bullet}} \ar[r] \ar[d, "\pi"] & \prod_{j \in J} \Higgs(C_j, \widetilde{\nu}^j_{\bullet}, \widetilde{\sigma}^j_{\bullet})^{nilp, \cO_{\bullet}} \ar[d, "\prod_{j} \vec{res}_j"] \\ \prod_{i \in I}(\mathfrak{gl}^{nilp}/\GL_N)  \ar[r, "\vec{\Delta}^-"] & \prod_{j \in J} \left(\mathfrak{gl}^{nilp}/\GL_N\right)^{n_j}
	\end{tikzcd}
\]
Since there are finitely many nilpotent conjugacy classes in $\mathfrak{gl}_N$, as we run over the set $\Sigma$ of all $I$-tuples $\cN_{\bullet}=(\cN_i)_{i \in I}$ of nilpotent conjugacy classes we obtain a stratification \[\Higgs(C, \sigma_{\bullet})^{nilp, \cO_{\bullet}}= \bigsqcup_{(\cN_i)_{i \in I} \in \Sigma} \pi^{-1}\left(\prod_{i \in I}(\cN_i/\GL_N)\right).\]
It suffices to show that each $\pi^{-1}\left(\prod_{i \in I}(\cN_i/\GL_N)\right)$ is equidimensional of expected dimension $N^2(g-1)+\frac{1}{2}\sum_{i=1}^n\dim(\cO_i)$. Note that $\pi^{-1}\left(\prod_{i \in I}(\cN_i/\GL_N)\right)$ fits into the following Cartesian diagram
 \begin{equation} \label{diagram 1}
	\begin{tikzcd}[ampersand replacement=\&]
		\pi^{-1}(\prod_{i \in I}(\cN_i/\GL_N)) \ar[r] \ar[d, "\pi"] \& \prod_{j \in J} \Higgs(C_j, \widetilde{\nu}^j_{\bullet}, \widetilde{\sigma}^j_{\bullet})^{nilp, \cO_{\bullet}} \ar[d, "\prod_{j} \vec{res}_j"] \\ \prod_{i \in I}(\cN_i/\GL_N)  \ar[r, "\vec{\Delta}^-"] \& \prod_{i \in I} \left(\mathfrak{gl}^{nilp}/\GL_N\right)^2
	\end{tikzcd}
\end{equation}
where we are using the identification $\prod_{i \in I} \left(\mathfrak{gl}^{nilp}/\GL_N\right)^2 \cong \prod_{j \in J} \left(\mathfrak{gl}^{nilp}/\GL_N\right)^{n_j}$. The fiber product 
 \[
	\begin{tikzcd}
		\prod_{j \in J} \Higgs(C_j, \widetilde{\nu}^j_{\bullet}, \widetilde{\sigma}^j_{\bullet})^{nilp, \cO_{\bullet}, \cN_{\bullet}} \ar[r] \ar[d] & \prod_{j \in J} \Higgs(C_j, \widetilde{\nu}^j_{\bullet}, \widetilde{\sigma}^j_{\bullet})^{nilp, \cO_{\bullet}} \ar[d, "\prod_{j} \vec{res}_j"] \\ \prod_{i \in I}\left(\cN_i/\GL_N\right)\times \left(-\cN_i/\GL_N\right)  \ar[r] & \prod_{i \in I} \left(\mathfrak{gl}^{nilp}/\GL_N\right)^2
	\end{tikzcd}
\]
parametrizes tuples of nilpotent meromorphic Higgs bundles $(E, \psi: E \to \otimes \omega_{C_j}^{log})$ in $\Higgs(C_j, \widetilde{\nu}^j_{\bullet}, \widetilde{\sigma}^j_{\bullet})^{nilp, \cO_{\bullet}}$ such that for each node $x_i\in C$ the residues at the two preimages $p_i,q_i$ are in $\cN_i$ and $-\cN_i$ respectively (in this case $\cN_i = -\cN_i$, but we won't use this). In view of Diagram \eqref{diagram 1}, the stack $\pi^{-1}\left(\prod_{i \in I}(\cN_i/\GL_N)\right)$ fits into a Cartesian diagram
 \begin{equation} \label{diagram 2}
	\begin{tikzcd}[ampersand replacement=\&]
		\pi^{-1}\left(\prod_{i \in I}(\cN_i/\GL_N)\right) \ar[r] \ar[d, "\pi"] \& \prod_{j \in J}\Higgs(C_j, \widetilde{\nu}^j_{\bullet}, \widetilde{\sigma}^j_{\bullet})^{nilp, \cO_{\bullet}, \cN_{\bullet}} \ar[d, "\prod_{j} \vec{res}_j"] \\ \prod_{i \in I}(\cN_i/\GL_N)  \ar[r, "\vec{\Delta}^-"] \& \prod_{i \in I}\left(\cN_i/\GL_N\right)\times \left(-\cN_i/\GL_N\right)
	\end{tikzcd}
\end{equation}
Since $\prod_{i \in I}\left(\cN_i/\GL_N\right)\times \left(-\cN_i/\GL_N\right)$ is the classifying stack of a group (because all the $\cN_i$'s are orbits), it follows that both of the following morphisms 
\[\prod_{j \in J}\Higgs(C_j, \sigma_{\bullet}^j, \tau_{\bullet}^j)^{nilp, \cO_{\bullet}, \cN_{\bullet}} \to \prod_{i \in I}\left(\cN_i/\GL_N\right)\times \left(-\cN_i/\GL_N\right)\]
\[\vec{\Delta}^-: \prod_{i \in I}(\cN_i/\GL_N) \to \prod_{i \in I}\left(\cN_i/\GL_N\right)\times \left(-\cN_i/\GL_N\right)\]
are automatically flat of pure relative dimension. It follows directly from definition that the relative dimension of $ \vec{\Delta}^-: \prod_{i \in I}(\cN_i/\GL_N) \to \prod_{i \in I}\left(\cN_i/\GL_N\right)\times \left(-\cN_i/\GL_N\right)$ is given by $\sum_{i \in I} (N^2-\dim(\cN_i))$. On the other hand, by the claim in the case of a smooth curve, we have that $\prod_{j \in J}\Higgs(C_j, \widetilde{\nu}^j_{\bullet}, \widetilde{\sigma}^j_{\bullet})^{nilp, \cO_{\bullet}, \cN_{\bullet}}$ is equidimensional of dimension $N^2(\sum_{j \in J} g_j - |J|) + \sum_{i \in I} \dim(\cN_i) + \frac{1}{2} \sum_{i=1}^N \dim(\cO_i)$. We conclude from the Diagram \eqref{diagram 2} with flat and equidimensional morphisms that the stack $\pi^{-1}\left(\prod_{i \in I}(\cN_i/\GL_N)\right)$ is equidimensional of dimension $N^2(\sum_{j \in J} g_j -|J| +|I| +1-1) + \frac{1}{2}\sum_{i=1}^n \dim(\cO_i)$. Since the genus $g$ of the curve $C$ is given by $g=\sum_{j \in J} g_j -|J| +|I| +1$, this is equal to the desired dimension $N^2(g-1) + \frac{1}{2}\sum_{i=1}^n \dim(\cO_i)$.
\end{proof}

\begin{prop} \label{prop: syntomicity of the possion Hitchin fibration}
    The morphism $H: \GHiggs_{g,n}^N \to \cA$ is syntomic of pure relative dimension $N^2(g-1) + \binom{N}{2}\cdot n$.
\end{prop}
\begin{proof}
    Using \Cref{prop: syntomicity of the residue morphism}, we see that the morphism $\GHiggs_{g,n}^N \to \overline{\cM}_{g,n}$ is syntomic of pure relative dimension $N^2(2g-2+n)$ when $n>0$ (resp. syntomic of pure relative dimension $N^2(2g-2)+1$ when $n=0$). (We need to distinguish two cases because the Fuchs stack has dimension $0$ when $n=0$, and has dimension $-1$ when $n\geq 1$.)
    By the dimension computations in the proof of \Cref{lemma: morphism from hitchin to bun is affine}, it follows (by separate computations for $n>0$ and $n=0$) that $\dim(\GHiggs_{g,n}^N) - \dim(\cA) = N^2(g-1) + \binom{N}{2}\cdot n$. Since $\GHiggs_{g,n}^N$ is syntomic and $\cA$ is smooth (\Cref{lemma: morphism from hitchin to bun is affine}), it suffices to show that the fibers of $H: \GHiggs_{g,n}^N \to \cA$ are equidimensional of the expected fiber dimension $N^2(g-1) + \binom{N}{2}\cdot n$. By using the $\mathbb{G}_m$-action that scales the meromorphic Higgs field (as in the previous proof of \Cref{prop: syntomicity of the hitchin morphism}), we reduce to showing that the preimage $H^{-1}(0_{\overline{\cM}_{g,n}})$ of the zero section $0_{\overline{\cM}_{g,n}}: \overline{\cM}_{g,n} \to \cA$ is equidimensional of the expected dimension $N^2(g-1) + \binom{N}{2}\cdot n$. But note that $H^{-1}(0_{\overline{\cM}_{g,n}})$ can also be described as the preimage of the zero section of the morphism $H: \GHiggs_{g,n}^{\cO_{\bullet}} \to \cA^{\cO_{\bullet}}$ where we set all conjugacy classes in the tuple $\cO_{\bullet}$ to be regular nilpotent. The desired dimension statement then follows directly from \Cref{prop: syntomicity of the hitchin morphism}.
\end{proof}
    
\end{subsection}

\end{section}
\begin{section}{Moduli spaces and Harder-Narasimhan stratifications}

Our goal in this section is to identify an open semistable locus $\GHiggs_{g,n}^{N,ss} \subset \GHiggs_{g,n}^N$ that admits a good moduli space in the sense of \cite{alper-good-moduli}, thus constructing a moduli space for meromorphic Gieseker Higgs bundles. We achieve this using the notion of $\Theta$-semistability with respect to a numerical invariant $\mu$ on the stack $\GHiggs_{g,n}^N$ (as in \cite{heinloth-hilbertmumford, halpernleistner2018structure}). To show the existence of the moduli space, we use the intrinsic criteria developed in \cite{alper2019existence} and \cite[\S5]{halpernleistner2018structure}. 
As a byproduct of our results in this section, we also construct a Harder-Narasimhan stratification of the stack $\GHiggs_{g,n}^N$ (it is what is known as a $\Theta$-stratification in \cite{halpernleistner2018structure}). This is used to conclude the properness of the Hitchin morphism, via the semistable reduction theorem in \cite{alper2019existence}.

\begin{subsection}{Review of numerical invariants and \texorpdfstring{$\Theta$}{Theta}-stratifications}

Our main strategy for the construction of the moduli space follows the criteria developed in \cite[\S 5]{halpernleistner2018structure} along with the technique of infinite-dimensional GIT introduced in \cite{torsion-freepaper, gauged_theta_stratifications}. 
In this subsection, we give a very minimalistic review of some of the relevant notions. We refer the reader to \cite[\S2.5]{torsion-freepaper} and \cite{gomez2023guide} for more detailed expositions on the contents of this subsection, and \cite{halpernleistner2018structure} for a complete treatment. 

    For the rest of this subsection, let $\mathscr{X}$ be an algebraic stack locally of finite type and with affine diagonal over $k$.

\subsubsection*{$\Theta$-stratifications}
\begin{defn}
        We denote by $\Theta$ the quotient stack $\mathbb{A}^1_{k}/\mathbb{G}_m = \Spec(k[t])/\mathbb{G}_m$, where $\mathbb{G}_m$ acts via its standard linear action which induces weight $-1$ on the coordinate function $t$. For any $k$-scheme $T$, we set $\Theta_T:= \Theta \times T$. If $T = \Spec(R)$ for some ring $R$, then we also use the notation $\Theta_R := \Theta_T$.
    \end{defn}
    There is an open immersion $1: \Spec(k) \hookrightarrow \Theta$ corresponding to the open orbit of $1 \in \mathbb{A}^1_k(k)$. There is also morphism $0: \Spec(k) \to \Theta$ 
    %which 
    with
    closed image corresponding to the origin $0 \in \mathbb{A}^1_k$. We also denote by $0: B \mathbb{G}_m \hookrightarrow \Theta$ the closed immersion of the residual gerbe at the image of $0$.

    \begin{defn} \label{defn: filtration point stack} Let $K \supset k$ be a field extension, and fix $p \in \mathscr{X}(K)$. A $\Theta$-filtration of $p$ consists of the data of a morphism $f: \Theta_K \to \mathscr{X}$ along with an isomorphism $ f(1) \xrightarrow{\sim} p$ for the image of the open point $1 \in \Theta_K(K)$. We denote by $\text{Filt}(p)$ the set of isomorphism classes of filtrations of $p$.
    \end{defn}

    We denote by $\text{Filt}(\mathscr{X}) = \Map(\Theta, \, \mathscr{X})$ the stack of filtrations of points in $\mathscr{X}$. By \cite[Prop. 1.1.2]{halpernleistner2018structure}, $\text{Filt}(\mathscr{X})$ is an algebraic stack locally of finite type and with affine diagonal over $k$. There is a morphism $\text{ev}_1: \text{Filt}(\mathscr{X}) \to \mathscr{X}$ given by evaluating at $1 \in \Theta(k)$.

 \begin{defn} \label{defn: theta stratification}
A $\Theta$-stratification of $\mathcal{X}$ consists of a collection of open substacks $\mathfrak{S}_c \subset \Filt(\mathscr{X})$ indexed by a totally ordered set $\Gamma$ such that the following are satisfied:

\begin{enumerate}[(1)]
\item For all $c$, the induced morphism $\ev_1: \mathfrak{S}_c \hookrightarrow \Filt(\mathscr{X}) \to \mathscr{X}$ is a locally closed immersion.
\item For any two distinct elements $c \neq c' \in \Gamma$, the locally closed substacks $\ev_1(\mathfrak{S}_c) \subset \mathscr{X}$ and $\ev_1(\mathfrak{S}_{c'}) \subset \mathscr{X}$ are disjoint. 
\item For any given $c \in \Gamma$, the union $\mathscr{X}_{\leq c} := \bigcup_{c' \leq c} \ev_1(\mathfrak{S}_{c'})$ is open in $\mathscr{X}$, and $\mathscr{X} = \bigcup_{c \in \Gamma} \mathscr{X}_{\leq c}$.

\item For every point $p \in \mathscr{X}$, the set $\left\{ c \in \Gamma \, \mid \, p \in \mathscr{X}_{\leq c}\right\}$ has a minimal element.
\end{enumerate}
\end{defn}

One may think of a $\Theta$-stratification as a stratification of the stack $\mathcal{X}$ where each locally closed stratum has a ``moduli-theoretic interpretation": it parametrizes pairs $(p,f)$ where $p$ is a point of $\mathscr{X}$ and $f$ is a $\Theta$-filtration of $p$ satisfying certain open conditions in the stack of filtrations $\text{Filt}(\mathscr{X})$. 

\begin{remark}
    The reader should be aware that the definition of $\Theta$-stratification that is usually found in the literature is made in terms of the open substacks $\mathscr{X}_{\leq c}$ without fixing a priori the choices of substacks of filtrations $\mathfrak{S}_c \subset \Filt(\mathscr{X})$. This does not make a big difference for practical purposes; we have decided to include the choice of $\mathfrak{S}_c \subset \Filt(\mathscr{X})$ for simplicity in our exposition.
\end{remark}

\begin{example}
    Let $C$ be a smooth projective connected curve over $k$, and let $N>0$. Consider the algebraic stack $\Bun_N(C)_0$ of rank $N$ and degree $0$ vector bundles on $C$. Given a $k$-point $p \in \Bun_N(C)(k)$ corresponding to a vector bundle $E$ on $C$, the Rees construction \cite[Prop. 1.0.1]{halpernleistner2018structure} induces a correspondence between the set of (isomorphism classes of) $\Theta$-filtrations of $p$ and the set of separated exhaustive decreasing $\mathbb{Z}$-indexed filtrations $(E_i)_{i \in \mathbb{Z}}$ by subbundles, i.e. the set of tuples $(E_i)_{i \in \mathbb{Z}}$ such that
    \begin{enumerate}[(1)]
        \item $E_i \subset E$ is a subbundle of $E$ (it is a subsheaf such that the quotient $E/E_i$ is a vector bundle).
        \item $E_i \supset E_{i+1}$ for all $i \in \mathbb{Z}$.
        \item $E_i =0$ for all sufficiently large $i$, and $E_i = E$ for all sufficiently small $i$.
    \end{enumerate}

    Hence, the stack $\Filt(\Bun_N(C)_0)$ parametrizes pairs $(E, E_{\bullet})$ of a degree $0$ rank $N$ vector bundle $E$ equipped with a filtration $E_{\bullet}$ satisfying conditions (1), (2) and (3) above. We say that the filtration $E_{\bullet}$ is a (weighted) Harder-Narasimhan filtration if for all $i \in \mathbb{Z}$ with $E_i / E_{i+1} \neq 0$, we have that $E_i / E_{i+1}$ is a semistable vector bundle and $i = \frac{\deg(E_i/E_{i+1})}{\text{rank}(E_i/E_{i+1})} \cdot N!$. In other words, the graded subquotients of the filtration are semistable, and the weights are proportional to the slopes of the subquotients. By openness of semistability and the local constancy of degrees and ranks of vector bundles, it follows that being a (weighted) Harder-Narasimhan filtration is an open condition in $\Filt(\mathscr{X})$. Hence, there is an open substack $\mathfrak{S} \subset \Filt(\Bun_N(C)_0)$ which parametrizes (weighted) Harder-Narasimhan filtrations. Given such a filtration $f= (E, E_{\bullet}) \in \mathfrak{S}(k)$, we define its degree of instability by $\mu(f) := \sum_{i \in \mathbb{Z}} i^2 \cdot \text{rank}(E_i/E_{i+1})$. Given a nonnegative integer $n\geq 0$, we denote by $\mathfrak{S}_{n} \subset \mathfrak{S}$ the open and closed substack of filtrations $f$ with $\mu(f)=n$. Then, the set of substacks $\mathfrak{S}_{n}\subset \Filt(\Bun_N(C)_0)$ indexed by the nonnegative integers $\Gamma= \mathbb{Z}_{\geq 0}$ is a $\Theta$-stratification of $\Bun_N(C)_0$.
\end{example}

\subsubsection*{Numerical invariants and semistability}

One of the main ways to construct a $\Theta$-stratification of a stack uses the notion of numerical invariant. Under certain hypotheses, such numerical invariant $\nu$ on the stack $\mathscr{X}$ will define a ``Harder-Narasimhan" $\Theta$-stratification, where we stratify based on the degree of instability of points of $\mathscr{X}$ with respect to $\nu$.

\begin{defn} \label{defn: general numerical invariant}
Let $\Gamma$ be an abelian group. A ($\Gamma$-valued) numerical invariant $\nu$ on the stack $\mathcal{X}$ is an assignment defined as follows. Let $K \supset k$ be a field extension and let $p \in \mathscr{X}(K)$. Let $\gamma: \left(\mathbb{G}_{m}^q\right)_K \rightarrow \text{Aut}(p)$ be a homomorphism of affine group schemes with finite kernel. 
The numerical invariant $\nu$ assigns to this data a scale-invariant function $\nu_{\gamma}: \mathbb{R}^q\setminus \{0\} \rightarrow \Gamma$ such that
\begin{enumerate}[(1)]
    \item $\nu_{\gamma}$ is unchanged under field extensions $K \subset K'$.
    \item $\nu$ is locally constant in algebraic families.
    \item Given a homomorphism $\phi: \left(\mathbb{G}^w_{m}\right)_K \rightarrow \left(\mathbb{G}^q_{m}\right)_K$ with finite kernel, the function $\nu_{\gamma \circ \phi}$ is the restriction of $\nu_{\gamma}$ along the inclusion $\mathbb{R}^w \cong X_*(\mathbb{G}_m^w)_{\mathbb{R}} \hookrightarrow X_*(\mathbb{G}_m^q)_{\mathbb{R}} \cong \mathbb{R}^q$ induced by $\phi$.
\end{enumerate}
\end{defn}

For any field extension $K \supset k$, we say that a $\Theta$-filtration $f: \Theta_K \to \mathscr{X}$ is nondegenerate if the restriction $f|_0: B\left(\mathbb{G}_m\right)_K \to \mathscr{X}$ induces a homomorphism of automorphism groups $\left(\mathbb{G}_m\right)_K \to \Aut(f(0))$ with finite kernel. 

For our purposes in this paper, numerical invariants will take values in the totally ordered group $\Gamma = \mathbb{R}[m, \epsilon]$ defined as follows.
\begin{notn}
    We denote by $\mathbb{R}[m, \epsilon]$ the ring of polynomials in two variables $m,\epsilon$ with coefficients in $\mathbb{R}$. We view this as a totally ordered abelian group under addition as follows:
  \begin{enumerate}[(1)]
      \item Given two polynomials $p,g \in \mathbb{R}[m]$, we say that $p \geq g$ if for all sufficiently large $m \gg0$ we have $p(m) \geq g(m)$.
      \item Given two polynomials $p,g \in \mathbb{R}[m, \epsilon]$, we say that $p \geq g$ if for all sufficiently small value $h >0$, we have an inequality $p(m,h) \geq g(m,h)$ in $\mathbb{R}[m]$. 
  \end{enumerate}
  For example, we have $m^2 - m> m^2 - m - \epsilon >0$ and $m - \epsilon m^2>0$, whereas $-m+\epsilon <0$.
\end{notn}

From now on, we implicitly assume that all numerical invariants are valued in the group $\Gamma = \mathbb{R}[m, \epsilon]$.

Given a numerical invariant $\nu$ on $\mathscr{X}$, we regard $\nu$ as an $\mathbb{R}[m, \epsilon]$-valued function on the set of nondegenerate filtrations by defining $\nu(f) \vcentcolon = \nu_{f|_{0}}(1) \in \bR[m, \epsilon]$.

\begin{defn} \label{defn: semistability numerical invariant}
    Let $\nu$ be a numerical invariant on $\mathcal{X}$. Let $K \supset k$ be an algebraically closed field extension. Let $p \in \mathscr{X}(K)$. We say that $p$ is semistable if all nondegenerate filtrations $f$ of $p$ satisfy $\nu(f) \leq 0$. Otherwise, we say that $p$ is unstable. We denote by $\mathscr{X}^{\nu \dash ss}$ the set of equivalence classes of semistable geometric points in the topological space $|\mathscr{X}|$.
\end{defn}

    \begin{defn}
        Let $p \in \mathscr{X}(K)$ be a geometric point. If $p$ is unstable, then a nondegenerate filtration $f$ of $p$ is called a Harder-Narasimhan filtration if $\nu(f) \geq \nu(f')$ for any other nondegenerate filtration $f'$ of $p$.

    By convention, if $p$ is semistable we say that its Harder-Narasimhan filtration is the trivial filtration $\Theta_K \to \Spec(K) \xrightarrow{p} \mathscr{X}$.
    \end{defn}

    We say that a numerical invariant $\nu$ on $\mathscr{X}$ defines a $\Theta$-stratification if the following conditions hold:
    \begin{enumerate}[(1)]
        \item every unstable geometric point $x \in \mathscr{X}(K)$ has a Harder-Narasimhan filtration that is unique up to pre-composing with a ramified covering $\Theta_K \xrightarrow{[n]} \Theta_K$ given by $t \mapsto t^n$.
        \item For each $c\in\mathbb{R}[m, \epsilon]_{\geq 0}$, the set of Harder-Narasimhan filtrations $f$ in $\Filt(\mathscr{X})$ with $\nu(f)=c$ is open in $\Filt(\mathscr{X})$.
        There is an open substack of such filtrations $\mathfrak{S}_c \subset \Filt(\mathscr{X})$ such that the collection $\left(\mathfrak{S}_c\right)_{c \in \mathbb{R}[m, \epsilon]}$ is a $\Theta$-stratification of $\mathscr{X}$.
    \end{enumerate}  

\subsubsection*{Numerical invariants via line bundles and norms on graded points}
In this subsection, we explain how to define numerical invariants from the data of polynomial families of line bundles and a norm on graded points of the stack.

\begin{defn}\label{defn: nondegenerate}
    Fix a positive integer $q \geq 1$ and a field extension $K \supset k$. A ($K$-valued) $q$-graded point of $\mathscr{X}$ is a morphism $B\left(\mathbb{G}^q_{m}\right)_K \to \mathscr{X}$. We say that $g: B\left(\mathbb{G}^q_{m}\right)_K \to \mathscr{X}$ is nondegenerate if the corresponding homomorphism $\gamma: \left(\mathbb{G}^q_{m}\right)_K \rightarrow \text{Aut}(g({\Spec(K))})$ has finite kernel.
\end{defn}

\begin{defn}[Rational quadratic norm on graded points]
\label{defn: rational_quadratic_norm}
A rational quadratic norm $b$ on graded points of $\mathscr{X}$ is an assignment defined as follows. Let $K \supset k$ be a field extension, and let $p \in \mathscr{X}(K)$. Let $\gamma: \left(\mathbb{G}^q_m\right)_K \to \text{Aut}(p)$ be a homomorphism with finite kernel. Then $b$ assigns to this datum a positive definite rational quadratic norm $b_{\gamma}: \mathbb{R}^q \to \mathbb{R}$ such that
\begin{enumerate}[(1)]
    \item $b_{\gamma}$ is unchanged under field extensions $K \subset K'$.
    \item $b$ is locally constant in algebraic families.
    \item Given a homomorphism $\phi: \left(\mathbb{G}^w_{m}\right)_K \rightarrow \left(\mathbb{G}^q_{m}\right)_K$ with finite kernel, $b_{\gamma \circ \phi}$ is the restriction of $b_{\gamma}$ along the inclusion $\mathbb{R}^w \cong X_*(\mathbb{G}_m^w)_{\mathbb{R}} \hookrightarrow X_*(\mathbb{G}_m^q)_{\mathbb{R}} \cong \mathbb{R}^q$ induced by $\phi$.
\end{enumerate}
\end{defn}

    Let $\cL$ be a line bundle on $\mathscr{X}$. For any morphism $g: B\left(\mathbb{G}^q_{m}\right)_K \rightarrow \mathscr{X}$, we may view the pullback line bundle $g^*(\cL)$ as a character of $\left(\mathbb{G}^q_{m}\right)_K$. Using the natural identification $X^*(\mathbb{G}^q_{m}) \cong \mathbb{Z}^q$, we may view $g^*(\cL)$ as a $q$-tuple of integers $(w^{(i)})_{i=1}^q$, which we call the weight of $g^*(\mathcal{L})$. 

    When $q=1$, we use the notation $\wt(g^*\cL) := w^{(1)}$ to denote the unique weight as above. If $f: \Theta_K \to \mathcal{X}$ is a $\Theta$-filtration, then set $\wt(\cL)(f) := \wt((f|_0)^*(\cL))$.

    \begin{defn} \label{defn: polynomial family of line bundles}
        A $\mathbb{Z}$-indexed tuple of line bundles $(\cL_m)_{m \in \mathbb{Z}}$ on $\mathscr{X}$ is called a polynomial family of line bundles if there is a polynomial $h(t) \in \Pic(\mathscr{X})_{\mathbb{Q}}[t]$ with coefficients in the rational Picard group $\Pic(\mathscr{X})_{\mathbb{Q}}$ such that for all $m \in \mathbb{Z}$ we have equality of $\mathbb{Q}$-line bundles $\cL_m = h(m)$.
    \end{defn}

    Given a polynomial family of line bundles $(\cL_m)_{m \in \mathbb{Z}}$ and a morphism $g: B\left(\mathbb{G}^q_{m}\right)_K \rightarrow \mathscr{X}$, we may extend the notion of weight by linearity and define a $q$-tuple $\wt(g^*(\cL_{m})) = (w^{(i)})_{i=1}^q \in \mathbb{Q}[m]^q$ of polynomials in $\mathbb{Q}[m]$.

      Similarly as above, when $q=1$ we use the notation $\wt(g^*\cL_m) := w^{(1)}$, and if $f: \Theta_K \to \mathscr{X}$ is a $\Theta$-filtration, then set $\wt(\cL_m)(f) := \wt((f|_0)^*(\cL_m))$.

\begin{defn}
    Given a polynomial family of line bundles $\mathcal{L}_m$ on $\mathcal{X}$, we define for every $g: B\left(\mathbb{G}^q_{m}\right)_K \rightarrow \mathscr{X}$ an $\mathbb{R}$-linear function $\wt(\cL_m)_g : \mathbb{R}^q \to \mathbb{R}[m]$ given by
\[ \wt(\cL)_g\left((r_i)_{i=1}^q\right) = \sum_{i=1}^q r_i \cdot w^{(i)} .\]
\end{defn}

\begin{defn} \label{nont: numerical invariant in terms of line bundle and norm on graded points}
    Given two polynomial families of line bundles $\cL_m$ and $\cG_m$ on $\mathscr{X}$ and a rational quadratic norm on graded points $b$, we denote by $\nu = (\wt(\cL_m) - \epsilon \wt(\cG_m))/\sqrt{b}$ the numerical invariant on $\mathscr{X}$ defined as follows. For all nondegenerate $g: B\left(\mathbb{G}^q_{m}\right)_K \to \mathscr{X}$ with corresponding morphism $\gamma: \left(\mathbb{G}^q_{m}\right)_K \to \text{Aut}(g(\Spec(k)))$, we set
\[ \nu_{\gamma}(\vec{r}) = \frac{-\wt(\cL_m)_{g}(\vec{r})- \epsilon \wt(\cG_m)_g(\vec{r})}{\sqrt{b_{\gamma}(\vec{r})}}.\]
\end{defn}

\begin{remark}
    In the case of the stack of vector bundles $\Bun(C)$ on a fixed smooth projective connected curve $C$, it is sufficient to consider a numerical invariant valued in the real numbers $\mathbb{R}$, which is defined similarly as above using a single line bundle $\cL$ (the determinant of cohomology line bundle). In our setup for this paper, we need two additional parameters $m$ and $\epsilon$ for our numerical invariant (see \Cref{NumInv}). The ``dominant term'' given by $-\wt(\cL_m)/\sqrt{b}$ in our case is a $\mathbb{R}[m]$-valued numerical invariant which measures Gieseker semistability of Higgs torsion-free sheaves. In order to obtain a reasonable notion of semistability, we include the additional term $- \epsilon \wt(\cG_m)/\sqrt{b}$ which can be thought of as a ``formal'' deformation of the stability condition coming from an ample polynomial family of line bundles on a relative stack of stable maps into a Quot scheme (cf. \cite{gauged-maps-projective}).
\end{remark}

\subsubsection*{Intrinsic criteria for moduli spaces and $\Theta$-stratifications}

Our main goal in this subsection is to explain the hypotheses in the following theorem.
\begin{thm}[{\cite[Theorem B]{halpernleistner2018structure}}] \label{thm: theta stability paper theorem}
Let $\nu = (-\wt(\cL_m) - \epsilon \wt(\cG_m))/\sqrt{b}$ be a numerical invariant on $\mathscr{X}$ defined by a pair of polynomial families of line bundles $\cL_m$ and  $\cG_m$, and a rational norm on graded points $b$ as in \Cref{nont: numerical invariant in terms of line bundle and norm on graded points}. Then, the following hold:
\begin{enumerate}[(1)]

 \item If $\nu$ is strictly $\Theta$-monotone (\Cref{defn: strictly theta monotone and STR monotone}), then it defines 
 a $\Theta$-stratification of $\mathscr{X}$ if and only if it satisfies the HN boundedness condition (\Cref{defn: HN boundedness}). In particular, the semistable locus $\mathscr{X}^{\nu\dash ss}$ is an open substack of $\mathscr{X}$.
    \item Assume that $\nu$ satisfies all the conditions in (1), so it defines a $\Theta$-stratification. In addition, suppose that $\nu$ is strictly $S$-monotone (\Cref{defn: strictly theta monotone and STR monotone}) and that the semistable locus $\mathscr{X}^{\nu \dash ss}$ can be written as a disjoint union of quasicompact open substacks. Then $\mathscr{X}^{\nu \dash ss}$ admits a separated good moduli space in the sense of \cite{alper-good-moduli}.
    \item Suppose that all the conditions in (2) are satisfied. Furthermore, suppose that $\mathscr{X}$ satisfies the existence part of the valuative criterion for properness for complete discrete valuation rings over $k$ \cite[\href{https://stacks.math.columbia.edu/tag/0CLK}{Tag 0CLK}]{stacks-project}. Then, every connected component of $\mathscr{X}^{\nu \dash ss}$ admits a good moduli space that is proper over $k$.
        
\end{enumerate}
\end{thm}

Let us explain the meaning of strict monotonicity. We first need to introduce some notation.
\begin{notn}
    Let $\kappa$ be a field and let $a \geq 1$ be an integer. We denote by $\mathbb{P}^1_{\kappa}[a]$ the $\mathbb{G}_m$-scheme $\mathbb{P}^{1}_{\kappa}$ equipped with the $\mathbb{G}_m$-action determined by the equation $t \cdot [x:y] = [t^{-a}x : y]$. We set $0 = [0: 1]$ and $\infty = [1:0]$.
\end{notn} 

\begin{notn}\label{notn: rigidified theta}
Let $R$ be a complete discrete valuation ring over $k$. We set $Y_{\Theta_{R}} \vcentcolon = \mathbb{A}^1_{R}$ equipped with its standard $\mathbb{G}_m$-action.  
\end{notn}
The scheme $Y_{\Theta_R}$ contains a unique $\mathbb{G}_m$-invariant closed point. We denote this point by $\mathfrak{o}$.

\begin{notn}\label{notn: rigidified str}
Let $R$ be a complete discrete valuation ring over $k$. Choose a uniformizer $\varpi$ of $R$. We set $Y_{\overline{ST}_R} := \Spec(R[s,t]/(st-\varpi))$ equipped with the $\bG_m$-action that assigns $t$ weight $-1$ and $s$ weight $1$.
\end{notn}
The isomorphism class of the $\mathbb{G}_m$-scheme $Y_{\overline{ST}_{R}}$ is independent of the choice of uniformizer $\varpi$. Note that $Y_{\overline{ST}_R}$ contains a unique $\bG_m$-invariant closed point cut out by the ideal $(s,t)$. We denote this point by $\mathfrak{o}$.

\begin{defn}[Strict monotonicity] \label{defn: strictly theta monotone and STR monotone} A numerical invariant $\nu$ on $\mathcal{X}$ is strictly $\Theta$-monotone (resp. strictly $S$-monotone) if the following condition holds. 

Let $R$ be any complete discrete valuation ring over $k$, and set $\mathscr{Y}= \Theta_R$ (resp. set $\mathscr{Y}:= \overline{ST}_R$). Choose a map $\varphi: (Y_\mathscr{Y} \setminus \mathfrak{o})/\, \mathbb{G}_m \rightarrow \mathscr{X}$. Then, after possibly replacing $R$ with a finite DVR extension, there exists a reduced and irreducible $\mathbb{G}_m$-equivariant Deligne-Mumford stack $\Sigma$ with maps $f: \Sigma \rightarrow Y_{\mathscr{Y}}$ and $\widetilde{\varphi}: \Sigma/ \mathbb{G}_m \rightarrow  \mathscr{X}$ such that the following hold:
\begin{enumerate}[({M}1)]
    \item $f$ is proper, $\mathbb{G}_m$-equivariant, and its restriction induces an isomorphism $f : \, \Sigma_{Y_{\mathscr{Y}} \setminus \mathfrak{o}} \xrightarrow{\sim} Y_{\mathscr{Y}} \setminus \mathfrak{o}$.
    \item The following diagram commutes
\begin{figure}[H]
\centering
\begin{tikzcd}
  \left(\Sigma_{Y_{\mathscr{Y}} \setminus \mathfrak{o}}\right)/ \, \mathbb{G}_m \ar[rd, "\widetilde{\varphi}"] \ar[d, "f"'] & \\   (Y_\mathscr{Y} \setminus \mathfrak{o})/ \, \mathbb{G}_m \ar[r, "\varphi"'] &  \mathscr{X}
\end{tikzcd}
\end{figure}
    \item Let $\kappa$ denote a finite extension of the residue field of $R$. For any $a \geq 1$ and any finite $\mathbb{G}_m$-equivariant morphism $\mathbb{P}^1_{\kappa}[a] \to \Sigma_{\mathfrak{o}}$, we have $\nu\left( \;\widetilde{\varphi}|_{\infty / \mathbb{G}_m} \;\right) >  \nu\left(\; \widetilde{\varphi}|_{0 / \mathbb{G}_m} \;\right)$.
\end{enumerate}
\end{defn}

Finally, the following is the definition of the HN boundedness condition in \Cref{thm: theta stability paper theorem}
\begin{defn}[HN boundedness]\label{defn: HN boundedness}
We say that a numerical invariant $\nu$ on $\mathscr{X}$ satisfies the HN boundedness condition if the following condition is always satisfied:

        Let $T$ be an affine Noetherian $k$-scheme, and choose a morphism $g: T \to \mathscr{X}$. Then, there exists a quasicompact open substack $\mathscr{U}_T \subset \mathscr{X}$ such that the following condition (*) holds.

            \noindent (*) For all geometric points $\overline{t} \to T$ with residue field $K \supset k$ and all nondegenerate $\Theta$-filtrations $f$ of the point $g(\overline{t})$ with $\nu(f)>0$, there exists another $\Theta$-filtration $f'$ of $g(\overline{t})$ with $f'(0) \in \mathscr{U}_T$ and $\nu(f')\geq \nu(f)$. 
\end{defn}

\end{subsection}
	\begin{subsection}{Definition of the numerical invariant}
	Let $(\cC, \sigma_{\bullet}) \to \overline{\cM}_{g,n}$ denote the universal stable $n$-pointed curve. Recall that the log canonical line bundle $ \omega_{\cC/\overline{\cM}_{g,n}}^{log}$ is relatively ample. This line bundle induces a notion of rank of torsion-free sheaves on fibers of $\mathcal{C}$, as in \cite[\S 2.3]{torsion-freepaper}. 
\begin{defn}[Rank of torsion-free sheaves]
    Let $(C, \sigma_{\bullet})$ be a stable curve defined over a field extension $K \supset k$. For any coherent sheaf $\cF$ on $C$, we denote by $\rk(\cF)$ the coefficient of degree 1 of the Hilbert polynomial of $\cF$ with respect to $\omega_{C/K}^{log}$. In particular, we have 
    \[\chi\left(\cF\otimes \left(\omega_{C/K}^{log}\right)^{\otimes m}\right) = \rk(\cF) m + \chi(\cF)
    \]
    for all $m$.
\end{defn}

For any field extension $K \supset k$ and any positive integer $h$, graded points $B(\mathbb{G}_m)^h_{K} \to \Coh^{tf}_{g,n}$ correspond to the data of an $n$-pointed stable curve $(C, \sigma_{\bullet})$ along with a $\bZ^h$-graded coherent sheaf $\bigoplus_{\vec{w} \in \bZ^h} \cF_{\vec{w}}$ on $C$ such that each graded piece $\cF_{\vec{w}}$ is torsion-free. Using this description, we define a rational quadratic norm on graded points as in \Cref{defn: rational_quadratic_norm} (see also \cite[Defn. 4.1.12]{halpernleistner2018structure}).
\begin{defn}[Quadratic form on the stack of
torsion-free sheaves]
We denote by $b$ the nondegenerate rational quadratic norm on graded points of $\Coh^{tf}_{g,n}$ defined as follows. For every field extension $K \supset k$, positive integer $h$, and nondegenerate point $g: B(\mathbb{G}_m)^h_{K} \to \Coh^{tf}_{g,n}$ corresponding to a $\mathbb{Z}^h$-graded sheaf $\bigoplus_{\vec{w} \in \bZ^h} \cF_{\vec{w}}$, we set
\[ b_g: \mathbb{R}^h \to \mathbb{R}, \; \; \; \; \vec{r} \mapsto \sum_{\vec{w} \in \bZ^h} \left(\vec{w} \cdot \vec{r} \right)^2  \rk(\F_{\vec{w}}). \]
where $\vec{w} \cdot \vec{r} = \sum_{j=1}^h w_j r_j$ denotes the standard inner product of the vectors $\vec{w} = (w_j)_{j=1}^h$ and $\vec{r} = (r_j)_{j=1}^h$.
\end{defn}

We also denote by $b$ the corresponding quadratic norm on $\GHiggs_{g,n}^N$ obtained by pulling back $b$ under the morphism $\GHiggs_{g,n}^N \to \GBun_{g,n}^N \to \Coh^{N}_{g,n} \subset \Coh^{tf}_{g,n}$. Since this composition morphism $\GHiggs_{g,n}^N \to \Coh^{tf}_{g,n}$ is representable (by \Cref{prop: gieseker vector bundles morphism proper} and \Cref{lemma: morphism from hitchin to bun is affine}), it follows that the quadratic norm $b$ is nondegenerate on $\GHiggs_{g,n}^N$.

For the following definition, we let $\pi_{\Coh}: \Coh^{tf}_{g,n} \times_{\overline{\cM}_{g,n}} \cC \to \Coh^{tf}_{g,n}$ be the first projection, and we let $\cF_{univ}$ denote the universal family of torsion free sheaves on $\Coh^{tf}_{g,n} \times_{\overline{\cM}_{g,n}} \cC$.
\begin{defn}[Determinant line bundle, cf. {\cite[Defn. 2.6]{torsion-freepaper}}] 
    For any $m \in \mathbb{Z}$, we denote by $\cL_{\det, m}$ the line bundle on $\Coh^{tf}_{g,n}$ given by
    \[ \cL_{\det, m} := \det\left(R(\pi_{\Coh})_*\left( \cF_{univ}\otimes (\omega_{\pi_{\Coh}}^{log})^{\otimes m}\right)\right) \]
\end{defn}

\begin{defn}[Gieseker line bundles, cf. {\cite[Defn. 2.12]{torsion-freepaper}}]
    We define the following line bundles on $\Coh^{tf}_{g,n}$:
    \begin{enumerate}[(1)]
        \item $b_1 := \cL_{\det,1} \otimes \cL_{\det,0}^{\vee}$.
        \item $\cL_{Gies, m} := \cL_{\det,m}^{\vee} \otimes b_1^{\otimes \overline{p}(m)}$,
where $\overline{p}(m) = m + \chi(\cF)/\rk(\cF)$ denotes the reduced Hilbert polynomial thought of as a locally constant function on $\Coh^{tf}_{g,n}$.
    \end{enumerate} 
    
    We shall also denote by $\cL_{Gies,m}$ pullback line bundles on $\GHiggs_{g,n}^N$ under the composition $\GHiggs_{g,n}^N \to \Higgs^{tf}_{g,n} \to \Coh^{tf}_{g,n}$.
\end{defn} 

\begin{remark}
    Note that $\cL_{Gies, m}$ is the dual $L_m^{\vee}$ of the line bundle in \cite[Defn. 2.12]{torsion-freepaper}).
\end{remark}

It was shown in \cite[Prop. 5.16]{torsion-freepaper} that $\Theta$-semistability for the family of line bundles $\cL_{Gies, m}$ recovers the notion of Gieseker semistability with respect to the polarization $\omega_{\mathcal{C}/\overline{\mathcal{M}}_{g,n}}^{log}$.

  \begin{defn}[Cornalba line bundles]
      We shall denote by $\cL_{Cor,m}$ the pullback of the corresponding line bundle on $\GBun_{g,n}^N$ (\Cref{defn: cornalba bundles}) under the forgetful morphism $\GHiggs^N_{g,n} \to \GBun_{g,n}^N$.
  \end{defn}
		
		We will need the following ampleness lemma.
		\begin{lemma} \label{lemma: ampleness of L_cor}
		Let $T$ be a quasicompact scheme equipped with a morphism $T \to \Higgs^{N}_{g,n}$. Then for all sufficiently large $m\gg0$ the pullback of the line bundle $\mathcal{L}_{Cor,m}$ is $T$-ample on the proper scheme $\GHiggs_{g,n}^{N} \times_{\Higgs^{N}_{g,n}} T$.
		\end{lemma}
		\begin{proof} It follows from \Cref{prop: gauged maps projective ampleness} that the pullback of $\mathcal{L}_{Cor,m}$ is $T$-ample on the scheme $\GBun_{g,n}^{N} \times_{\Coh^{N}_{g,n}} T$ for all sufficiently large $m \gg 0$. The lemma follows, because $\GHiggs_{g,n}^{N} \times_{\Higgs^{N}_{g,n}} T$ is a closed subscheme of $\GBun_{g,n}^{N} \times_{\Coh^{N}_{g,n}} T$ by the proof of \Cref{lemma: properness of gieseker hitchin stack over torsion-free hitchin stack}
		\end{proof}

  Now that we have polynomial families of line bundles $\cL_{Cor,m}$ and $\cL_{Gies,m}$ and a nondegenerate rational quadratic norm on our stack, we may define a numerical invariant as explained in \Cref{nont: numerical invariant in terms of line bundle and norm on graded points}.
		
		\begin{defn}\label{NumInv}
			We define the ($\mathbb{R}[m, \epsilon]$-valued) numerical invariant $\mu$ on $\GHiggs_{g,n}^{N}$ as 
			\[\mu := \frac{-\wt(\mathcal{L}_{Gies,m}) - \epsilon \cdot \wt(\cL_{Cor,m})}{\sqrt{b}}.\]	
		\end{defn}
	\end{subsection}
	\begin{subsection}{Proof of monotonicity}
		 We proceed to prove monotonicity for the numerical invariant $\mu$ on $\GHiggs_{g,n}$. The argument uses the strategy of infinite-dimensional GIT as in \cite{torsion-freepaper, gauged_theta_stratifications, gauged-maps-projective}.

  \begin{prop} \label{prop: monotonicity of the numerical invariant}
			The numerical invariant $\mu$ 
   given in \Cref{NumInv}
   is strictly $\Theta$-monotone and strictly $S$-monotone on the stack $\GHiggs_{g,n}$ (as in \Cref{defn: strictly theta monotone and STR monotone}).
		\end{prop}
		\begin{proof}
    		Let $R$ be a complete discrete valuation ring over $k$. Set $Y$ to be one of the $\mathbb{G}_m$-schemes $Y_{\Theta_R}$ or $Y_{\overline{ST}_R}$ as defined in \Cref{notn: rigidified theta} and \Cref{notn: rigidified str}, and we denote by $\mathfrak{o}$ the unique $\mathbb{G}_m$-invariant closed point in $Y$. Suppose that we are given a $\mathbb{G}_m$-equivariant $k$-morphism $Y \setminus \mathfrak{o} \to \GHiggs_{g,n}^N$. Consider the composition $Y \setminus \mathfrak{o} \to \GHiggs_{g,n}^N \to \Higgs_{g,n}^{tf}$, which classifies a pair $(C, \sigma_{\bullet}, \mathcal{F}, \psi)$ consisting of a $\mathbb{G}_m$-equivariant $Y \setminus \mathfrak{o}$-family of stable $n$-pointed curves and $\mathbb{G}_m$-equivariant family $(\mathcal{F}, \psi)$ of torsion-free meromorphic Higgs sheaves on $C$. Since the stack $\overline{\cM}_{g,n}$ is proper Deligne-Mumford, we can uniquely complete $(C, \sigma_{\bullet})$ to an equivariant family of stable $n$-pointed curves $(\widetilde{C}, \widetilde{\sigma}_{\bullet}) \to Y$.
			
			The proof of \cite[Thm. 4.10]{torsion-freepaper} shows that there exists an ind-projective ind-scheme $\Gr \to Y$ equipped with an action of $\mathbb{G}_m$ and a $\mathbb{G}_m$-equivariant morphism $\Gr \to \Higgs^{tf}(\widetilde{C}/Y)$ to the $Y$-stack $\Higgs^{tf}(\widetilde{C}/Y)$ of torsion-free meromorphic Higgs sheaves on the fibers of the family of pointed curves $(\widetilde{C}, \widetilde{\sigma}_{\bullet}) \to Y$. Furthermore, this morphism fits into the following $\mathbb{G}_m$-equivariant commutative diagram
			\begin{figure}[H]
				\centering
				\begin{tikzcd}
					Y \setminus \mathfrak{o} \ar[r] \ar[dr] &  \Gr \ar[d] \ar[r, "\xi"]& \Higgs^{tf}(\widetilde{C}/Y) \ar[d] \ar[r] & \Higgs^{tf}_{g,n} \\   & Y \ar[r, symbol = \xrightarrow{\sim}] & Y &
				\end{tikzcd}
			\end{figure}
\noindent such that the composition $Y \setminus \mathfrak{o} \to \Higgs^{tf}_{g,n}$ yields $(C, \sigma_{\bullet}, \mathcal{F}, \psi)$. By construction, any two field-valued points $p_1,p_2 \in \Gr$ mapping to the same fiber $y \in Y$ yield torsion-free meromorphic Higgs sheaves $\xi(p_1), \xi(p_2)$ on $\widetilde{C}$ that are generically isomorphic. By \cite[Cor. 3.28]{torsion-freepaper}, the line bundle $\cL_{Gies,m}|_{\Gr}$ is ample on each of the $Y$-projective strata of $\Gr$ for $m \gg0$. Since $Y \setminus \mathfrak{o}$ is quasicompact, the morphism $Y \setminus \mathfrak{o} \to \Gr$ factors through one of these $Y$-projective strata $\Gr^{i}$. Consider the fiber product $\Gr \times_{\Higgs_{g,n}^{tf}}$, which is a $\Gr^i$-proper Deligne-Mumford stack equipped with an action of $\mathbb{G}_m$ and fitting into the following $\mathbb{G}_m$-equivariant commutative diagram

		\begin{figure}[H]
			\centering
			\begin{tikzcd}
				& \Gr^i \times_{\Higgs^{tf}_{g,n}} \GHiggs_{g,n}^N \ar[r] \ar[d] & \GHiggs_{g,n}^N \ar[d] \\
				Y \setminus \mathfrak{o} \ar[r] \ar[dr]  \ar[ur] &  \Gr^i \ar[d] \ar[r] & \Higgs_{g,n}^{tf} \\   & Y &
			\end{tikzcd}
		\end{figure}
		Here the composition $Y \setminus \mathfrak{o} \to \Higgs^{tf}_{g,n}$ recovers the original $\mathbb{G}_m$-equivariant family we started with. For all sufficiently large $m\gg0$, the pullback of the formal line bundle $\cL_{Gies,m} + \epsilon \cL_{Cor,m}$ on the $Y$-proper Deligne-Mumford stack $\Gr^i \times_{\Higgs^{tf}_{g,n}} \GHiggs_{g,n}^N$ is eventually (for $m\gg0$) positive relative to $Y$ in the sense of \cite[Defn. 2.12]{gauged_theta_stratifications}. This is because $\cL_{Gies,m}$ is eventually $Y$-ample on $\Gr^i$, and $\mathcal{L}_{Cor,m}$ is eventually $\Gr^i$-ample on $\Gr^i \times_{\Higgs^{tf}_{g,n}} \GHiggs_{g,n}^N$ by \Cref{lemma: ampleness of L_cor}.
		
		Let $\Sigma \hookrightarrow \Gr^i \times_{\Higgs_{g,n}^{tf}} \GHiggs_{g,n}^N$ denote the scheme theoretic closure of $Y \setminus \mathfrak{o}$. By construction, $\Sigma$ is a $Y$-proper Deligne-Mumford stack equipped with a compatible $\mathbb{G}_m$ action, and the pullback of the formal line bundle $\Sigma|_{\cL_{Gies,m} + \epsilon \cL_{Cor,m}}$ is eventually (for $m\gg0$) positive relative to $Y$. Also, since the fiber $\Sigma_\mathfrak{o}$ maps to $\mathfrak{o}\in Y$, by the properties of $\Gr$ described above, we see that any two $\mathbb{G}_m$-equivariant points of $\Sigma_\mathfrak{o}$ yield graded torsion-free meromorphic Higgs sheaves that are generically isomorphic on $\widetilde{C}_\mathfrak{o}$. Since the norm $b$ only depends on the ranks of the graded components, we see that any two $\mathbb{G}_m$-equivariant points $p_1,p_2: \Spec(K) \to \Sigma_\mathfrak{o} \to \GHiggs_{g,n}^N$ yield graded points $g_1, g_2$ of $\GHiggs_{g,n}^N$ with the same norm $b(g_1) = b(g_2)$. Using this and the positivity of the formal polynomial line bundle $\cL_{Gies,m} + \epsilon \cL_{Cor,m}$, one concludes the proof of strict monotonicity in the same way as in the proof of \cite[Prop. 3.6]{gauged_theta_stratifications} (or in the end of the proof of \cite[Thm. 4.9]{torsion-freepaper}).
		\end{proof}
	\end{subsection}

 \begin{subsection}{Filtrations and semistability}
Recall the definition of a $\Theta$-filtration of a point in a stack given in \Cref{defn: filtration point stack}. It turns out that we can describe $\Theta$-filtrations in $\GHiggs_{g,n}^N$ via the morphism $\varpi: \GHiggs_{g,n}^{N} \to \Higgs_{g,n}^N$ given by pushing forward to the stabilization.
		\begin{lemma} \label{lemma: comparison of filtrations}
			Let $p$ be a geometric point of $\GHiggs_{g,n}^{N}$. Then the morphism $\varpi$ establishes a bijection between $\Theta$-filtrations of $p$ and $\Theta$-filtrations of the image $\varpi(p)$.
		\end{lemma}
		\begin{proof}
			This follows from \cite[Thm. B.1]{gauged_theta_stratifications}, because the morphism $\varpi: \GHiggs_{g,n}^{N} \to \Higgs^{N}_{g,n}$ is schematic and proper (\Cref{lemma: properness of gieseker hitchin stack over torsion-free hitchin stack}).
		\end{proof}

     Fortunately, $\Theta$-filtrations in the stack $\Higgs^{N}_{g,n}$ can be described explicitly via the Rees construction \cite[Prop. 1.0.1]{halpernleistner2018structure} (see also \cite[\S4.1]{torsion-freepaper}). Let $K \supset k$ be a field extension, and fix a $K$-point $p=(\overline{C}, \overline{\sigma}_{\bullet}), \cF,\psi)$ of $\Higgs^{N}_{g,n}$. A $\Theta$-filtration of $p$ corresponds to a $\mathbb{Z}$-indexed sequence $(\cF_l)_{l \in \mathbb{Z}}$ of subsheaves $\cF_l \subset \cF$ satisfying the following properties:
     \begin{enumerate}[(1)]
         \item $\cF_{l+1} \subset \cF_l$.
         \item $\cF_l=0$ for $l\gg0$ sufficiently large, and $\cF_l = \cF$ for $l\ll 0$ sufficiently small.
         \item $\cF_l/\cF_{l+1}$ is a pure sheaf of dimension $1$ (torsion-free).
         \item $\psi(\cF_l) \subset \cF_l \otimes \omega_{\overline{C}/K}^{log}$.
     \end{enumerate}
    Associated to every nondegenerate $\Theta$-filtration $f: \Theta_K \to \GHiggs_{g,n}^N$, we have a corresponding polynomial
    \[ \mu(f) = \frac{-\wt(\cL_{Gies,m})(f) - \epsilon \cdot \wt(\cL_{Cor,m})(f)}{\sqrt{b(f|_0)}} \in \mathbb{Q}[m,\epsilon].\]

The $\epsilon$-leading term $-\wt(\cL_{Gies,m})/\sqrt{b}$ is the pullback under $\varpi$ of the Gieseker numerical invariant for $\Higgs^{N}_{g,n}$ defined in \cite{torsion-freepaper}. The computations in \cite[\S4.1]{torsion-freepaper} allow us to give an explicit description of this invariant. Let $p = (\overline{C}, \sigma_{\bullet}, \cF, \psi)$ be a geometric point of $\Higgs^N_{g,n}$, and let $f$ be a $\Theta$-filtration of $p$ corresponding to a sequence of subsheaves $(\cF_l)_{l \in \mathbb{Z}}$. Then, we have
\[ -\frac{\wt(\cL_{Gies,m})(f)}{\sqrt{b(f|_0)}} = \frac{\sum_{l \in \mathbb{Z}} l \cdot \rk(\cF_l/\cF_{l+1})\cdot \left( \frac{\chi(\cF_l/\cF_{l+1})}{\rk(\cF_l/\cF_{l+1})} - \frac{\chi(\cF)}{\rk(\cF)}\right)}{\sqrt{\sum_{l \in \mathbb{Z}} l^2 \cdot \rk(\cF_l/\cF_{l+1})}}.\]

\begin{remark}
    Note in particular that $-\wt(\cL_{Gies,m})(f)/\sqrt{b(f|_0)}$ is a constant (as a polynomial in $m$).
\end{remark}

\begin{notn}[Gieseker (semi)stability of torsion-free meromorphic Higgs sheaves]
    Let $(\overline{C}, \sigma_{\bullet}, \cF, \psi)$ be a torsion-free meromorphic Higgs sheaf. Recall that $(\cF, \psi)$ is called Gieseker semistable with respect to the polarization $\omega_{\overline{C}}^{log}$
(in the sense of polynomial semistability in \cite[\S3]{Simpson-repnI} when viewed a $\Lambda$-module) if for all subsheaves $\cG \subset F$ with $\psi(\cG) \subset \cG \otimes \omega_{\overline{C}}^{log}$, we have $\chi(\cG)/\rk(\cG) \leq \chi(\cF)/\rk(\cF)$. 

We say that $(\cF, \psi)$ is Gieseker stable if for any subsheaf $\cG \subset \cF$ such that $\psi(\cG) \subset \cG\otimes \omega_{\overline{C}}^{log}$ and $\cF/\cG$ is a nonzero torsion-free sheaf, we have strict inequality $\chi(\cG)/\rk(\cG) < \chi(\cF)/\rk(\cF)$.
\end{notn}
     
    We end this subsection with the following observation. Recall the definition of semistability with respect to a numerical invariant in \Cref{defn: semistability numerical invariant}.
     \begin{coroll} \label{coroll: semistable implies semistable torsion-free pushforward}
			Let $(C,\sigma_{\bullet}, E, \psi)$ be a geometric point of $\GHiggs_{g,n}^{N}$ that is semistable with respect to the numerical invariant $\mu_{Gies}$. Then the pushforward $(\varphi_*(E), \varpi_*(\psi)) \in \Higgs^{N}_{g,n}$ is a Gieseker semistable meromorphic Higgs sheaf on the stabilization $(\overline{C}, \overline{\sigma}_{\bullet})$.
		\end{coroll}
		\begin{proof}
			This is an immediate consequence of \Cref{lemma: comparison of filtrations}, because the $\epsilon$-leading term $-\wt(\cL_{Gies,m})/\sqrt{b}$ of the numerical invariant $\mu$ is the pullback of the numerical invariant on $\Higgs^{tf}_{g,n}$ that recovers Gieseker semistability (cf. \cite[\S5]{torsion-freepaper}).
		\end{proof}

 \end{subsection}
 
	\begin{subsection}{HN boundedness and boundedness of the semistable locus}
		
		\begin{prop} \label{prop: hn boundedness}
			The numerical invariant $\mu$ on $\GHiggs_{g,n}^{N}$ satisfies the HN boundedness condition (see \Cref{defn: HN boundedness}).
		\end{prop}
		\begin{proof}
			Let $T$ be a Noetherian scheme and $g:T \to \GHiggs_{g,n}^{N}$ be a family parametrized by $T$. Without loss of generality, we may assume that there is some $\chi \in \mathbb{Z}$ such that the image $g(T)$ lands in the open and closed substack $\GHiggs_{g,n}^{N,\chi} \subset \GHiggs_{g,n}^N$ that parametrizes meromorphic Gieseker Higgs bundles of whose underlying vector bundle has Euler characteristic $\chi$. The composition $T \to \GHiggs_{g,n}^{N} \xrightarrow{\varpi} \Higgs_{g,n}^{N}$ lands inside the analogous open and closed substack $\Higgs_{g,n}^{N,\chi}$. Furthermore, the image of $T$ belongs to an open stratum $\left(\Higgs_{g,n}^{N, \chi}\right)^{\leq \gamma}$ for some $\gamma \in \mathbb{R}$, where we are considering the strata with respect to the numerical invariant $-\wt(\cL_{Gies,m})/\sqrt{b}$ of torsion-free meromorphic Higgs sheaves (viewed as $\Lambda$-modules) considered in \cite{torsion-freepaper}. Recall that $\left(\Higgs_{g,n}^{N, \chi}\right)^{\leq \gamma}$ consists of those torsion-free meromorphic Higgs sheaves of fixed Hilbert polynomial (determined by $N,\chi$) such that the supremum of the numerical invariant over all $\Theta$-filtrations is smaller than or equal to $\gamma$. This imposes an upper bound on the reduced Hilbert polynomial of any meromorphic Higgs subsheaf (viewed as a $\Lambda$-submodule). Hence, by \cite[Lemma 3.3]{Simpson-repnI} combined with \cite[Thm. 4.4]{langer-boundedness}, it follows that $\left(\Higgs_{g,n}^{N, \chi}\right)^{\leq \gamma}$ is of finite type. By the quasicompactness of $\varpi$ (\Cref{lemma: properness of gieseker hitchin stack over torsion-free hitchin stack}), we conclude that $\left(\Higgs_{g,n}^{N,\chi}\right)^{\leq \gamma}$ is of finite type over $k$.
   
   To finish the proof, it suffices to show that the open quasicompact substack $\mathscr{U}_T := \varpi^{-1}\left(\left(\Higgs_{g,n}^{N,\chi}\right)^{\leq \gamma}\right)$ satisfies the condition (*) in the definition of HN boundedness. This amounts to showing that, in order to maximize the numerical invariant $\mu$ of filtrations of an unstable geometric point $p: \overline{t} \to T \to \GHiggs_{g,n}^N$, it suffices to consider filtrations $f: \Theta_{\overline{t}} \to \GHiggs_{g,n}^{N}$ such that $\varpi \circ f$ lands in $\left(\Higgs_{g,n}^{N, \chi}\right)^{\leq \gamma}$. We prove this by considering two separate cases:

 (a) The torsion-free meromorphic Higgs sheaf $\varpi(p)$ is Gieseker unstable. Since the Gieseker numerical invariant $-\wt(\cL_{Gies,m})/\sqrt{b}$ in \cite{torsion-freepaper} defines a $\Theta$-stratification on $\Higgs_{g,n}^{N,\chi}$, there is a unique maximal filtration up to scaling $f$ of $\varpi(p)$ that has maximal numerical invariant. Since $\varpi$ induces a bijection between filtrations of $p$ and $\varpi(p)$ (up to scaling) and since the $\epsilon$-leading term of $\mu$ is the pullback $-\wt(\cL_{Gies,m})/\sqrt{b}$ of the Gieseker numerical invariant on $\Higgs_{g,n}^N$, it follows that the unique inverse image $\widetilde{f}_{max} = ``\varpi^{-1}(f)"$ must be the unique maximizer (up to scaling) of $\mu$ among filtrations of $p$. Since the maximizing filtration $f$ lands in the same stratum as $\varpi(p)$ in the $\Theta$-stratification of $\Higgs_{g,n}^{N,\chi}$, we conclude that indeed $\varpi(\widetilde{f}_{max}) = f$ lands in $\left(\Higgs_{g,n}^{N,\chi}\right)^{\leq \gamma}$.
				
(b) The torsion-free meromorphic Higgs sheaf $\varpi(p)$ is Gieseker semistable. This means that for all filtrations $f$ of $p$, the $\epsilon$-leading term of $-\mu(\cL_{Gies,m})/\sqrt{b}$ will be nonpositive. Therefore, in order to maximize the numerical invariant, we can restrict to filtrations $f$ such that the $\epsilon$-leading term of $-\mu(\cL_{Gies,m})/\sqrt{b}$ is $0$. This will happen if and only if the composition $\varpi \circ f$ lands in the semistable locus $\left(\Higgs_{g,n}^{N,\chi}\right)^{\leq 0} \subset \left(\Higgs_{g,n}^{N,\chi}\right)^{\leq \gamma}$, as desired.
		\end{proof}
		Since $\mu$ is strictly $\Theta$-monotone (Proposition \ref{prop: monotonicity of the numerical invariant}) and satisfies HN boundedness, we conclude by \Cref{thm: theta stability paper theorem}(1) that $\mu$ defines a $\Theta$-stratification on $\GHiggs_{g,n}^{N}$. In particular, the set of $\mu$-semistable geometric points form an open substack of $\GHiggs_{g,n}^{N}$, which we will denote by $\GHiggs_{g,n}^{N,ss}$. 
  
    For each $\chi \in \mathbb{Z}$, we denote by $\GHiggs^{N,\chi}_{g,n} \subset \GHiggs_{g,n}^N$ the open and closed substack parametrizing $T$-families $(C, \sigma_{\bullet}, E, \psi)$ such that all the $T$-fibers of the family of vector bundles $E$ have Euler characteristic $\chi$. We denote by $\GHiggs_{g,n}^{N,\chi,ss} \subset \GHiggs_{g,n}^{N,ss}$ the corresponding open and closed substack of $\GHiggs_{g,n}^{N,\chi,ss}$ where the Euler characteristic of the underlying vector bundle is $\chi$.
		
		\begin{prop} \label{prop: boundedness of semistable locus}
			For all $\chi \in \mathbb{Z}$, the stack $\GHiggs_{g,n}^{N,\chi,ss}$ is of finite type over $k$.
		\end{prop}
		\begin{proof}
		Let $\Higgs^{N,\chi, ss}_{g,n} \subset \Higgs^{tf}_{g,n}$ denote the open substack of Gieseker semistable torsion-free meromorphic Higgs pairs of generic rank $N$ and Euler characteristic $\chi$. As we saw over the course of the proof of \Cref{prop: hn boundedness}, the stack $\Higgs^{N,\chi, ss}_{g,n}$ is of finite type over $k$. By the quasicompactness of the morphism $\varpi$ (\Cref{lemma: properness of gieseker hitchin stack over torsion-free hitchin stack}), the preimage $\varpi^{-1}\left(\Higgs^{N,\chi, ss}_{g,n}\right)$ is also of finite type over $k$. By \Cref{coroll: semistable implies semistable torsion-free pushforward}, the semistable locus $\GHiggs_{g,n}^{N,\chi,ss}$ is an open substack of $\varpi^{-1}\left(\Higgs^{N,\chi, ss}_{g,n}\right)$, and hence it is of finite type over $k$.
		\end{proof}
	\end{subsection}
	\begin{subsection}{Moduli spaces of meromorphic Gieseker Higgs bundles}
		\begin{thm} \label{thm: theta-stratification and moduli space for logarithmic gieseker higgs}
       The following statements hold:
			\begin{enumerate}[(i)]
				\item The numerical invariant $\mu$ given in \Cref{NumInv}
    defines a $\Theta$-stratification on the stack $\GHiggs_{g,n}^{N}$.
				\item Fix an integer $\chi \in \mathbb{Z}$. Then the substack $\GHiggs_{g,n}^{N,\chi,ss}$ admits a relative good moduli space $\MGHiggs_{g,n}^{N,\chi}$ over $\overline{\cM}_{g,n}$ that is separated and of finite type over $\overline{\cM}_{g,n}$. Moreover, the induced Hitchin morphism $H: \MGHiggs_{g,n}^{N,\chi} \to \cA$ is flat and proper.
			\end{enumerate}
		\end{thm}
		\begin{proof}
		In \Cref{prop: monotonicity of the numerical invariant} we have seen that $\mu$ is strictly $\Theta$-monotone and strictly $S$-monotone. Propositions \ref{prop: hn boundedness} and \ref{prop: boundedness of semistable locus} show that $\mu$ satisfies the HN boundedness condition, and that the closed and open substack $\GHiggs_{g,n}^{N,\chi,ss}$ of $\GHiggs_{g,n}^{N,ss}$ is of finite type. Moreover, by \Cref{prop: properness of the hitchin morphism}, the Hitchin morphism $H: \GHiggs_{g,n}^{N} \to \cA$ satisfies the existence part of the valuative criterion for properness. Therefore, everything in the theorem except for the flatness of $H$ follows from \Cref{thm: theta stability paper theorem}, which we may apply \'etale-locally over the base $\overline{\cM}_{g,n}$. The flatness of $H: \MGHiggs_{g,n}^{N,\chi} \to \cA$ follows from the flatness of the corresponding morphism from the stack $H: \GHiggs_{g,n}^N \to \cA$ (\Cref{prop: syntomicity of the possion Hitchin fibration}) combined with \cite[Thm. 4.16(ix)]{alper-good-moduli}.
		\end{proof}
	\end{subsection}

 As an immediate consequence, we get the following theorem.
 \begin{thm} \label{thm: moduli space for closure of stack with fixed residues}
     Fix an integer $\chi \in \mathbb{Z}$ and a tuple $\cO_{\bullet}$ of nilpotent conjugacy classes as in Subsection \ref{subsection: hitchin morphism for fixed nilpotent residues}. Then the semistable locus $\overline{\GHiggs}_{g,n}^{\cO_{\bullet}, \chi, ss} \subset \overline{\GHiggs}_{g,n}^{\cO_{\bullet}, \chi}$ is an open substack which admits a relative good moduli space $\overline{\MGHiggs}^{\cO_{\bullet}, \chi}_{g,n}$ over $\overline{\cM}_{g,n}$. There is an induced morphism $H:\overline{\MGHiggs}^{\cO_{\bullet},\chi}_{g,n} \to \cA^{\cO_{\bullet}}$ that is proper and flat.
 \end{thm}
 \begin{proof}
     The openness of semistability and the existence of the relative good moduli space are direct consequences of \Cref{thm: theta-stratification and moduli space for logarithmic gieseker higgs}, because $\overline{\GHiggs}_{g,n}^{\cO_{\bullet}} \hookrightarrow \GHiggs_{g,n}^N$ is a closed substack. Furthermore, we have an induced closed immersion $\overline{\MGHiggs}^{\cO_{\bullet}, \chi}_{g,n} \subset \MGHiggs^{N,\chi}_{g,n}$. It follows from \Cref{prop: factorization of the hitchin morphism with fixed nilpotent residues} that the morphism $\overline{\MGHiggs}^{\cO_{\bullet}, \chi}_{g,n} \hookrightarrow \MGHiggs^{N,\chi}_{g,n} \to \cA$ factors uniquely through $\cA^{\cO_{\bullet}} \to \cA$. Since $\overline{\MGHiggs}^{\cO_{\bullet}, \chi}_{g,n} \hookrightarrow \MGHiggs^{N,\chi}_{g,n} \to \cA$ is proper, it follows that $H: \overline{\MGHiggs}^{\cO_{\bullet},\chi}_{g,n} \to \cA^{\cO_{\bullet}}$ is proper. Moreover, since the morphism $\overline{\GHiggs}_{g,n}^{\cO_{\bullet}} \to \cA^{\cO_{\bullet}}$ is flat (\Cref{prop: syntomicity of the hitchin morphism}),
     it follows by \cite[Thm. 4.16(ix)]{alper-good-moduli} that the induced morphism $H: \overline{\MGHiggs}^{\cO_{\bullet},\chi}_{g,n} \to \cA^{\cO_{\bullet}}$ is flat.
 \end{proof}

For the next discussion, we will need a notion of stable points. Notice that there is a canonical copy of $\mathbb{G}_m$ inside the group of automorphisms of any meromorphic Gieseker Higgs bundle $p=(C, \sigma_{\bullet}, E, \psi)$, namely, the group of automorphisms that acts by scaling $E$. In particular, we have $\dim(\Aut(p)) \geq 1$. By upper-semicontinuity of fiber dimension for the inertia $I_{\GHiggs_{g,n}^{N,\chi}} \to \GHiggs_{g,n}^{N,\chi}$, the locus $\mathscr{Z} \subset |\GHiggs_{g,n}^{N,\chi}|$ of points $p$ in $\GHiggs_{g,n}^{N,\chi}$ with $\dim(Aut(p))>1$ is closed. Consider the good moduli space morphism $\pi: \GHiggs_{g,n}^{N, \chi} \to \MGHiggs_{g,n}^{N,\chi}$. By \cite[Thm. 4.16(ii)]{alper-good-moduli}, the image $\pi(\mathscr{Z})$ is closed in $\MGHiggs_{g,n}^{N,\chi}$.

\begin{notn}[Moduli space of stable meromorphic Gieseker Higgs bundles]
    We set $\MGHiggs_{g,n}^{N,\chi,s}$ to be the open complement of $\pi(\mathscr{Z}) \subset \MGHiggs_{g,n}^{N,\chi}$.
\end{notn}

\begin{defn}[Stack of stable meromorphic Gieseker Higgs bundles]
    We denote by $\GHiggs_{g,n}^{N,\chi,s}$ the open complement of $\pi^{-1}(\pi(\mathscr{Z})) \subset \GHiggs_{g,n}^{N,\chi}$ as defined above.
\end{defn}

\begin{remark}
    Let $\GHiggs^{N,\chi,s}_{g,n} \subset \GHiggs_{g,n}^{N,\chi}$ denote the open substack of Gieseker stable torsion-free meromorphic Higgs sheaves. Since the group scheme of automorphisms of any Gieseker stable object in $\Higgs^{N,\chi,s}_{g,n}$ relative to $\overline{\cM}_{g,n}$ is isomorphic to $\mathbb{G}_m$, it follows that $\varpi^{-1}(\Higgs^{N,\chi,s}_{g,n}) \subset \GHiggs_{g,n}^{N,\chi,s}$, where $\varpi: \GHiggs_{g,n}^N \to \Higgs^{N}_{g,n}$ is the proper schematic morphism defined in \Cref{lemma: properness of gieseker hitchin stack over torsion-free hitchin stack}.
\end{remark}

\begin{prop} \label{prop: tame moduli space for stable Gieseker Higgs}
    The morphism $\pi: \GHiggs_{g,n}^{N,\chi,s} \to \MGHiggs_{g,n}^{N,\chi,s}$ is a tame moduli space morphism in the sense of \cite[\S7]{alper-good-moduli}.
\end{prop}
\begin{proof}
   It follows from construction that $\GHiggs_{g,n}^{N,\chi,s} \subset \GHiggs_{g,n}^{N,\chi}$ is a saturated open substack with respect to $\pi$ \cite[Defn. 6.1]{alper-good-moduli}, and hence $\pi: \GHiggs_{g,n}^{N,\chi,s} \to \MGHiggs_{g,n}^{N,\chi,s}$ is a good moduli space morphism. Since the automorphism groups of all points in $\GHiggs_{g,n}^{N,\chi,s}$ are $1$-dimensional, it follows from \cite[Prop. 9.1]{alper-good-moduli} that the geometric fibers of $\pi: \GHiggs_{g,n}^{N,\chi,s} \to \MGHiggs_{g,n}^{N,\chi,s}$ consist of a single point.
\end{proof}

\begin{defn}[Stable meromorphic Gieseker Higgs bundles with fixed residues]
    Fix an $n$-tuple of conjugacy classes $\cO_{\bullet}$. We set $\GHiggs_{g,n}^{\cO_{\bullet},\chi,s}:= \GHiggs_{g,n}^{\cO_{\bullet}} \cap \GHiggs_{g,n}^{N,\chi,s}$, where the intersection takes place inside $\GHiggs_{g,n}^{N,\chi}$.
\end{defn}

 \begin{thm}
     Fix an integer $\chi \in \mathbb{Z}$ and an $n$-tuple of conjugacy classes $\cO_{\bullet}$ as in Subsection \ref{subsection: hitchin morphism for fixed nilpotent residues}. The open substack of stable points $\GHiggs_{g,n}^{\cO_{\bullet},\chi,s} \subset \GHiggs_{g,n}^{\cO_{\bullet},\chi}$ admits a relative tame moduli space $\MGHiggs^{\cO_{\bullet}, \chi,s}_{g,n}$ over $\overline{\cM}_{g,n}$. There is an induced morphism $H: \MGHiggs^{\cO_{\bullet}, \chi,s}_{g,n} \to \cA^{\cO_{\bullet}}$ that is flat of pure relative dimension $N^2(2g-2)+2 + \sum_{i=1}^N\dim(\cO_i)$.
 \end{thm}
 \begin{proof}
     Let $\pi: \overline{\GHiggs}_{g,n}^{\cO_{\bullet},\chi} \to \overline{\MGHiggs}_{g,n}^{\cO_{\bullet},\chi}$ denote the good moduli space morphism. It follows from \Cref{prop: tame moduli space for stable Gieseker Higgs} that the open substack $\GHiggs_{g,n}^{\cO_{\bullet},\chi,s} \subset \overline{\GHiggs}_{g,n}^{\cO_{\bullet},\chi}$ is saturated with respect to $\pi$, and that the geometric fibers of $\pi: \GHiggs_{g,n}^{\cO_{\bullet},\chi,s} \to \overline{\MGHiggs}^{\cO_{\bullet},\chi}_{g,n}$ are singletons. Therefore, $\GHiggs_{g,n}^{\cO_{\bullet},\chi,s}$ admits the relative moduli space $\MGHiggs_{g,n}^{\cO_{\bullet},\chi,s}:= \pi(\GHiggs_{g,n}^{\cO_{\bullet},\chi,s})$, which is an open inside $\overline{\MGHiggs}_{g,n}^{\cO_{\bullet},\chi}$. The flatness of the morphism $H: \MGHiggs_{g,n}^{\cO_{\bullet},\chi,s} \to \cA^{\cO_{\bullet}}$ follows immediately from \Cref{thm: moduli space for closure of stack with fixed residues}. The dimension computation for the fibers follows from the dimension of the fibers of $H: \GHiggs_{g,n}^{\cO_{\bullet}} \to \cA^{\cO_{\bullet}}$ (\Cref{prop: syntomicity of the hitchin morphism}) and the fact that all the points of $\GHiggs_{g,n}^{\cO_{\bullet}}$ have $1$-dimensional automorphism groups.
 \end{proof}

 \end{section}

\footnotesize{\bibliography{nodal_hitchin.bib}}
\bibliographystyle{alpha}
\vspace{0.5cm}
  \noindent \textsc{Department of Mathematics, University of Pennsylvania,
209 South 33rd Street,
Philadelphia, PA 19104, USA}\par\nopagebreak
\noindent \textit{E-mail address}: \texttt{donagi@math.upenn.edu}\par\nopagebreak
  \noindent \textit{E-mail address}: \texttt{andresfh@sas.upenn.edu}
  
\end{document}